\documentclass[12pt,a4paper]{article}

\usepackage{a4,enumerate}
\usepackage{amsthm}
\usepackage{latexsym}
\usepackage{times}
\usepackage{pstricks,pst-coil,pst-plot}
\usepackage{eufrak}

\newdimen\templaenge

\DeclareMathAlphabet{\doba}{U}{msb}{m}{n}
\gdef\mC{\doba{C}}

\gdef\mN{\doba{N}}

\gdef\mR{\doba{R}}

\gdef\mT{\doba{T}}
\gdef\mZ{\doba{Z}}

\def\qed{{\leavevmode\unskip\nobreak\hfil\penalty 50\hskip 1em%
  \hbox{}\nobreak\hfil\lower 1pt\hbox{$\Box$\kern-.5pt}\parfillskip 0pt
  \finalhyphendemerits 0\par\bigbreak}}
\def\qedmath#1{\setbox0\hbox{$\displaystyle #1$}\templaenge=\textwidth\advance\templaenge by -\wd0%
\setbox1\hbox{$\Box$}\advance\templaenge by -2\wd1%
$$#1\hbox to0pt{\kern.5\templaenge$\Box$\kern-.5pt\hss}$$\par\bigbreak}

\def\al{{\alpha}}
\def\be{{\beta}}
\def\de{{\delta}}
\def\De{{\Delta}}
\def\om{{\omega}}

\def\la{{\lambda}}
\def\ka{{\kappa}}

\def\si{{\sigma}}
\def\Si{{\Sigma}}
\def\ga{{\gamma}}
\def\ep{{\varepsilon}}
\def\Ga{{\Gamma}}

\def\th{{\vartheta}}
\def\Th{{\Theta}}
\def\ph{{\varphi}}
\def\phi{{\varphi}}

\def\rh{{\rho}}

\def\na{{\nabla}}
\def\pa{{\partial}}
\def\el{{\ell}}

\def\cF{\mathcal{F}}

\def\cK{\mathcal{K}}
\def\cL{\mathcal{L}}
\def\cM{\mathcal{M}}

\def\cS{\mathcal{S}}
\def\cW{\mathcal{W}}
\def\cV{\mathcal{V}}

\def\ohne{-}
\def\ti{\tilde}
\def\witi{\widetilde}
\def\wihat{\widehat}

\def\lan{\langle}
\def\ran{\rangle}

\def\torus{T^2}
\def\Gd{\Ga_{xy}^*}
\def\will{\cW}

\let\mo=\mathopen
\let\mc=\mathclose

\def\ie{i.\thinspace e.\ \ignorespaces} 
\def\eg{e.\thinspace g.\ \ignorespaces}

\def\nummerarray#1#2{\par\noindent\setbox0\hbox{\rm (#1)}\setbox1\hbox{$#2$}\unhcopy0%
\dimen0=.5\textwidth \advance\dimen0 by -\wd0 \advance\dimen0 by -.5\wd1 \kern\dimen0 \unhcopy1}

\def\SO{\mathop{{\rm SO}}}
\def\SU{\mathop{{\rm SU}}}
\def\Spin{\mathop{{\rm Spin}}}
\def\GLp{{{\mathop{{\rm GL}}}^+}}
\def\GLpt{{\witi\GLp}}
\def\GLpte{{\witi{\mathop{{\rm GL}}}^+_1}}
\def\spinset{\mathop{{\mathfrak{Spin}}}}
\def\ker{\mathop{{\rm ker}}}
\def\grad{{\mathop{{\rm grad}}}}
\def\Span{\mathop{{\rm span}}}

\def\dim{\mathop{{\rm dim}}}
\def\vol{{\mathop{{\rm vol}}}}
\def\dvol{{\mathop{{\rm dvol}}}}
\def\area{{\mathop{{\rm area}}}}

\def\min{\mathop{{\rm min}}}
\def\max{\mathop{{\rm max}}}
\def\osc{\mathop{{\rm osc}}}

\def\id{\mathop{{\rm id}}}

\def\End{{\mathop{{\rm End}}}}

\def\Hom{{\mathop{{\rm Hom}}}}
\def\clott{C_{\textrm{\small Lott}}}
\def\cammann{C_{\textrm{\small Ammann}}}

\def\gflach{{g_0}}
\def\gbel{g}

\def\abel{{{\rm Area}_g}}
\def\afl{{{\rm Area}_0}}
\def\platz{{\,\mathord{.}\,}}

\def\res#1#2{{#1}\lower .11ex\hbox{$|$}\lower .644ex\hbox{$\scriptstyle #2$}}

\def\length{{\mbox{\rm length}}}

\def\Lnorm#1#2#3{\left\|#1\right\|_{L^{#2}(#3)}}
\def\lnorm#1#2{\left\|#1\right\|_{L^{#2}}}
\def\sys{{\mathop{\rm sys}}_{1}}
\def\spinsys{{\mathop{\textrm{\upshape spin-sys}}}_{1}}
\def\nssys{{\mathop{\textrm{\upshape nonspin-sys}}}_{1}}
\def\umax{{\max u}}
\def\umin{{\min u}}

\def\spinmod{\cM^{\textrm{\upshape spin}}}
\def\const{{\textit{const}}}

\long\def\komment#1{}

\def\proof#1{{\par\medbreak\noindent {\bf Proof\setbox0\hbox{#1}%
\ifdim\wd0=0pt .\else\ \ignorespaces #1.\fi}\enspace}}
\def\iop#1{{\par\medbreak\noindent {\bf Idea of proof\setbox0\hbox{#1}%
\ifdim\wd0=0pt .\else\ \ignorespaces #1.\fi}\enspace}}
\def\examples{{\noindent {\bf Examples. }\par\kern-\baselineskip}}

\newtheoremstyle{remarks}{3pt}{3pt}{}{}{\bfseries}{}{ }{}

\newtheorem{theorem}{\bf T{\footnotesize HEOREM}}[section]
\newtheorem{proposition}[theorem]{\bf P{\footnotesize ROPOSITION}}
\newtheorem{lemma}[theorem]{\bf L{\footnotesize EMMA}}
\newtheorem{corollary}[theorem]{\bf C{\footnotesize OROLLARY}}

\newtheorem*{diracproposition}{\bf P{\footnotesize ROPOSITION}~\ref{diracvergl}}
\newtheorem*{laplaceproposition}{\bf P{\footnotesize ROPOSITION}~\ref{laplacevergl}}

\theoremstyle{definition}
\newtheorem*{remark}{Remark}

\newtheorem*{definition}{Definition}
\newtheorem*{example}{Example}

\theoremstyle{remarks}

\def\eref#1{{\rm (\ref{#1})}}

\begin{document}
\title{Spectral estimates on $2$-tori}
\author{Bernd Ammann\footnote{ammann@math.uni-hamburg.de}}
\date{March 2000}
\maketitle

\begin{abstract}
We prove upper and lower bounds for the eigenvalues of the Dirac operator and the
Laplace operator on 2-dimensional tori. In particluar we give a lower bound
for the first eigenvalue of the Dirac operator for non-trivial spin structures.
It is the only explicit estimate for eigenvalues of the Dirac operator known so far
that uses information about the spin structure. 

As a corollary we obtain lower bounds for the Willmore functional 
of a torus embedded into $S^3$. 

In the final section we compare Dirac spectra for two different spin structures
on an arbitrary Riemannian spin manifold. 
\end{abstract}


{\bf Mathematics Classification 2000:}
Primary: 53C27, Secondary: 58J50 53C80

\section{Introduction}

The Dirac operator is an elliptic differential operator of order one playing an important role both in modern
physics and in mathematics. In physics, particles with non-integer spin, so-called fermions,
are described by the Dirac equation. Let us assume that the space-time $M$ is stationary, $M=\mR\times N$
and that the spatial component $N$ is compact and admits a spin structure.
Then stationary fermions have a wave function of the form
  $$ \Psi(t,x)=e^{iEt}\Psi_0(x)\qquad t\in\mR,\,x\in N$$
where $\Psi_0$ is an eigenspinor of $D_N^2$, the square of the Dirac operator on $N$, that belongs to the
eigenvalue~$\la$ . The energy $E$ and the eigenvalue $\la$ are related via the formula
  $$E^2=\la+m^2$$
with $m$ being the rest mass of the particle.
Knowing the spectrum therefore means knowing possible energies.
The first eigenvalue is of particular interest as it characterizes the energy of the state of lowest energy
--- the vacuum.
On an arbitrary Riemannian manifold, exact calculation of the spectrum is impossible, thus one tries try
to find bounds for the eigenvalues.

Bounding eigenvalues of the Dirac operator on a compact Riemannian manifold~$N$ 
is also an important tool in
differential geometry and topology. If $N$ is spin and carries a metric whose scalar curvature
is greater than or equal to $s_0>0$ at every point, then with
the help of the Schr\"odinger-Lichnerowicz formula
it is easy to prove that the first eigenvalue $\la_1$ of $D^2$ is bounded from below by $s_0/4$.
On the other hand Atiyah-Singer index theorem tells us
that positivity of the first eigenvalue of $D^2$ on a compact Riemannian
manifold~$N$ implies that the $\hat A$-genus vanishes.
Therefore any compact spin manifold admitting a positive scalar curvature metric has 
vanishing $\hat A$-genus.

Lower bounds for Dirac eigenvalues can also be applied to problems in classical
differential geometry. For any immersion
$F:N\to \mR^n$ of a compact manifold $N$, Christian B\"ar \cite{baer:98}
proved 
\begin{equation}\label{hzweiabschallg}
  \int_N |H|^2\geq \mu_1\area(N).
\end{equation}
Here $N$ carries the induced metric, $\mu_1$ is the first eigenvalue 
of the square of a twisted Dirac operator and $H$ is the mean curvature vector field 
of $F(N)\subset \mR^n$.
If $N$ is the 2-dimensional torus $\torus$, then the left hand side of \eref{hzweiabschallg} 
is the so-called Willmore functional. The Willmore conjecture  
states 
  $$\int_{\torus} |H|^2\geq 2\pi^2$$
for any immersion $F:\torus\to \mR^n$. This conjecture first appeared in \cite{willmore:65} for the case $n=3$.  
In the meantime the conjecture has been verified for several classes of immersions, for example for immersions with
rotational symmetry \cite{langer.singer:84} or for non-injective immersions \cite{li.yau:82}.
Nevertheless the conjecture remains open until now.
For further information on this conjecture the reader may read the introductions of 
\cite{topping:p98} or \cite{ammann:00}. 

Now assume for simplicity
that $F$ is an embedding and $F(\torus)\subset S^3\subset \mR^4$. 
In this case, the twisting bundle is trivial, and $\mu_1$ is the first eigenvalue of the square
of the classical Dirac operator associated to a non-trivial spin structure. 
Our goal is to use inequality~\eref{hzweiabschallg} in order to derive lower bounds 
for the Willmore
functional. If the induced metric on $\torus$ is flat, the spectrum of $D$ has been explicitely calculated
\cite{friedrich:84} and we obtain a lower bound for $\int_{\torus} |H|^2$.

Obtaining lower eigenvalue estimates for non-flat tori is much harder.
John Lott \cite[Proposition 1]{lott:86} proved the existence
of a constant $\clott>0$ depending on the spin-conformal type of the torus such that
\begin{equation}\label{lott}
  \mu_1\area\geq \clott.
\end{equation}
Unfortunately, Lott's article does not give an explicit value and it seems hard 
to express such a constant
$\clott$ in terms of meaningful geometric data. Lott's estimate uses the $L^p$-boundedness of zero order
pseudo-differential operators and Sobolev embedding theorems, 
hence corresponding constants are hard to interpret without using explicit
coordinates. 

The starting point of the author's PhD thesis \cite{ammanndiss} and of the present article is to 
find an explicit lower bound for $\mu_1$ that uses information about the spin structure.
All explicit lower estimates known before did not use any information
about the spin structure. 
\komment{
Even on the 2-torus there is no explicit lower eigenvalue bound depending on the spin structure. 
And because the 2-torus has $\mu_1=0$ for the trivial spin structure 
there was no non-trivial lower estimate
for $\mu_1$ on the $2$-dimensional torus for non-trivial spin structures.
}

For general compact Riemannian manifolds the problem of finding such estimates is rather difficult.
It is not clear 
at all what kind of data from the spin structure could be used in order 
to get an additional term in a lower eigenvalue
estimate.
Take for example a compact manifold with non-vanishing $\hat A$-genus. 
It has $\mu_1=0$ for any spin structure, thus the contribution of the spin structure
in the estimate has to vanish.

As the general case is hard to handle, most of the article will specialize to the 2-di\-men\-sio\-nal torus $\torus$.
By the uniformization theorem any 2-di\-men\-sio\-nal torus is conformally equivalent to a flat torus. 
We  use this fact in order to control the geometry. An important, but also very technical step 
for this is the estimate
of the oscillation of the conformal factor (Section~\ref{streckabschsec}).
Although our main goal was to find lower estimates for the Dirac eigenvalues, it turns out that 
this method gives upper and lower bounds for all eigenvalues both of the Laplace operator and the 
Dirac operator and for any spin structure. We prove different versions of the estimates. 
Theorem~\ref{spicotheo} for example states for the first eigenvalue $\mu_1$ of the square of the Dirac operator
\begin{equation}\label{meineabsch}
    \mu_1\area \geq \cammann \cdot \ka
\end{equation}
where $\cammann>0$ is an explicit constant depending on the spin-conformal class 
and $\ka\leq 1$ is a curvature expression that satisfies $\ka=1$ if the metric is flat.
This estimate is sharp for any flat metric.

In view of Lott's result \eref{lott}, it is tempting to conjecture that we can drop the curvature term, \ie 
$\mu_i\area \geq C_{\textrm{Ammann}}.$
This is false however: we can prove by example 
at the end of section~\ref{willmoreappl}
that for many spin-conformal structures the optimal constant 
in Lott's estimate is not attained by a flat torus.    

In section~\ref{willmoreappl} we will prove some lower bounds for the Willmore functional that are
strongly related to our lower estimates of the Dirac eigenvalues. 
In particular we prove for embeddings $\torus\to S^3$ that under a curvature
condition the Willmore functional converges to $\infty$ if the spin-conformal 
type of the embedding converges to one end of the spin-conformal moduli space (Corollary~\ref{willmoreinfty}).

The results in this paper about the Willmore conjecture are strongly related to another preprint of the author
\cite{ammann:00}. The results of the present article are stronger near one of the ends of the 
spin-conformal moduli space but they have other drawbacks. Namely, they do not generalize easily to higher codimensions
and they impose a restriction on the spin-conformal class.

The structure of the article is as follows: In section~\ref{statement} we will state our spectral estimates on 2-tori.
Sections~\ref{overview} to \ref{lastproofsection} provide proofs of
the statements in section~\ref{statement}.
We then apply Theorem~\ref{streckabsch} once again and derive an application 
to the Willmore functional 
that is related to our lower eigenvalue estimates.

Finally in section~\ref{beldim} we will prove a result for arbitrary spin manifolds $M$. 
Let $M$ carry two different spin structures $\th$ and $\th'$. 
The difference of these spin structures $\chi:=\th-\th'$ is an element in 
$H^1(M,\mZ_2)=\Hom_\mZ(H_1(M,\mZ),\mZ_2)$. Assume that $\chi$ vanishes on the torsion part of $H_1(M,\mZ)$.
We will define a norm $\lnorm{\chi}\infty$, the stable norm of $\chi$. 
We prove that the eigenvalues $(\rh_i)_{i\in\mZ}$ of the Dirac operator corresponding
to $\th$ and the eigenvalues $(\rh_i')_{i\in\mZ}$ corresponding to $\th'$ can be numbered so that
  $$|\rh_i-\rh'_i|\leq 2\pi\lnorm{\th-\th'}\infty.$$
If the spectrum is known for $\th$ and 
if  $|\rh_i|> 2\pi \lnorm{\th-\th'}\infty$ for any $i\in\mZ$,
then this yields a lower bound for any $\rh'_i$.

At the end of the introduction we want to mention some other publications 
that treat the interplay between spin structures and the spectrum of the Dirac operator.
However, they do not derive explicit eigenvalue bounds for generic metrics.
We will restrict to the most recent ones. For further references and a good overview of the subject
we refer to \cite{baer:p00}.

Dahl \cite{dahl:p99} shows that the difference of the eta-invariants corresponding to two 
different spin structures is an integer, 
if the difference of the spin structures viewed as an element in 
$\Hom_\mZ(H_1(M,\mZ),\mZ_2)$ vanishes on the torsion part.
B\"ar \cite{baer:p98} calculated the essential spectrum of 
hyperbolic 2- and 3-manifolds of finite volume. 
In these examples, the essential spectrum depends on the spin structure at the cusps.
Pf\"affle \cite{pfaeffle:p99} calculated the spectrum and the $\eta$-invariants of flat Bieberbach manifolds. 
These spectra also depend on the spin structure.

Several results in the present article already appeared in the author's PhD thesis \cite{ammanndiss}.

\section{Main results}\label{statement}

In this section we summarize our results about the spectra of Dirac and Laplace operators on 2-tori. 

The spectrum of the Dirac operator depends on the spin structure. 
At first, we recall some important facts about spin structures and introduce some notation.
Spin structures will be discussed in more detail in section~\ref{spinsection}.  

Let $M$ be a compact orientable manifold with vanishing second Stiefel-Whitney class $w_2(TM)=0$.
Such manifolds admit a spin structure. However, the spin structure is not unique in general.
The group $H^1(M,\mZ_2)$ acts freely and transitively on the set of spin structures $\spinset(M)$, \ie
$\spinset(M)$ is an affine space associated to the vector space $H^1(M,\mZ_2)$.
After fixing a spin structure and a Riemannian metric on $M$ we can define 
the spinor bundle $\Si M\to M$ and a Dirac operator 
$D:\Ga(\Si M)\to \Ga(\Si M)$.

We are mainly interested in the case $M=\torus$. 
The $2$-dimensional torus $\torus$ is spin. Because of $\#\spinset(\torus)=\# H^1(\torus,\mZ_2)=4$ 
there are 4 spin structures on $\torus$.
There is exactly one spin structure in $\spinset(\torus)$ for which $0$ lies in the spectrum of $D$, 
regardless of the underlying metric $g$. 
This spin structure will be called \emph{trivial} (see section~\ref{spinsection} for other characterizations).
We will identify the trivial spin structure with $0\in H^1(\torus,\mZ_2)$. 
This identification yields an identification of the affine space $\spinset(\torus)$ with $H^1(\torus,\mZ_2)$.
On the other hand, we will identify $H^1(\torus,\mZ_2)$ with $\Hom_\mZ(H_1(\torus,\mZ),\mZ_2)$. 
Hence spin structures on $\torus$ are in a canonical one-to-one relation to such homomorphisms. 
Frequently, we will use the term ``\emph{spin homomorphism}'' instead of ``spin structure'' in order to indicate
that we regard the spin structure as an element in $\Hom_\mZ(H_1(\torus,\mZ),\mZ_2)$.

If the torus $\torus$ carries a flat metric, it is very helpful to 
write the torus as $\mR^2/\Ga$ with a lattice $\Ga\cong H_1(\torus,\mZ)$. 
We always assume that $\mR^2/\Ga$ carries the metric induced by the Euclidean metric on $\mR^2$.
Let $\Ga^*$ be the lattice dual to $\Ga$.
Elements $\chi\in\Hom_\mZ(\Ga,\mZ_2)$ are represented by vectors $\al\in (1/2)\Ga^*$ with the property
  $$\chi(x) = (-1)^{2\al(x)}\qquad \forall x\in\Ga.$$
Note that $\chi$ determines $\al$ only up to elements in $\Ga^*$.

We define the function 
$\cS:[0,4\pi\mc[\times [0,\infty\mo[\times \mo]1,\infty\mc[\times\mo]0,\infty\mc]\to\mo]0,\infty]$ by 
  $$\cS(\cK,\cK',p,\cV):={p\over p-1}\,\left[ {\cK'\over 4\pi}+{1\over 2}\bigg|\log \left(1-{\cK\over 4\pi}\right)\bigg|+
    {\cK\over 8\pi -2\cK}\,\log\left({2\cK'\over \cK}\right)\right] + {\cK\cV\over 8}$$
for $\cK> 0$ and $\cS(0,\cK',p,\cV):=0$.

Let $\abel$ be the area of $(\torus,\gbel)$.

\begin{theorem}\label{maineins}
Let $(\torus,\gbel)$ be a Riemannian 2-torus with spin homomorphism $\chi$.
Choose a lattice $\Ga$ in $\mR^2$ with $\vol(\mR^2/\Ga)=1$ together with a conformal map $A:\mR^2/\Ga\to (\torus,\gbel)$.
Assume that $A^*(\chi)$ is represented by $\al\in(1/2)\Ga^*$.
Let $0\leq \el_0 \leq \el_1 \leq \el_2\leq\dots$ be the sequence of lengths of $\Ga^*+\al$ (with multiplicities), and let 
$(\mu_i\,|\,i=1,2,\dots)$ be the spectrum of $D^2$ on $(\torus,\gbel,\chi)$.

Then 
  $$e^{-2\osc u}\, 4\pi^2\, \el_{\left[{i-1\over 2}\right]}^2 \leq \mu_i \, \abel \leq e^{2\osc u}\, 
    4\pi^2\, \el_{\left[{i-1\over 2}\right]}^2.$$
If $\Lnorm{K_\gbel}1{\torus,\gbel}<4\pi$, then
\begin{equation}\label{hjkl} 
  \osc u \leq \cS\left(\Lnorm{K_\gbel}1{\torus,\gbel},\Lnorm{K_\gbel}p{\torus,\gbel}
           \abel^{1-(1/p)},p,\si_1(\torus,\gbel)^{-2}\right)
\end{equation}
with $\si_1(\torus,\gbel):=\inf\Big\{\length(\be)\,|\,\be\in\Ga\ohne\{0\}\Big\}$.
\end{theorem}

The number $\si_1(\torus,\gbel)$ is a conformal invariant of $(\torus,\gbel)$ which 
will be called \emph{cosystole}.

The most difficult step in the proof of this theorem is to find the estimate \eref{hjkl}. 
This step will be performed in 
Theorem~\ref{streckabsch}. For proving the above theorem, 
we will use the explicit formula for the spectra of flat tori 
(Proposition~\ref{diracflatspec}, \cite{friedrich:84}).
Another important tool for the proof is the following proposition. 

\begin{diracproposition}
Let $M$ be a compact manifold with two conformal metrics $\ti g$ and $g=e^{2u}\ti g$.
Let $D$ and $\witi D$ be the corresponding Dirac operators with respect to a common
spin structure.
We denote the eigenvalues of $D^2$ by
$\mu_1\leq\mu_2\leq\dots$ and the ones of ${\witi D}^2$
by $\witi\mu_1\leq\witi\mu_2\leq\dots$.

Then
  $$ \mu_i\,\min_{m\in M}e^{2u(m)}
     \leq \witi \mu_i \leq \mu_i\,  \max_{m\in M} e^{2u(m)} \qquad \forall i=1,2,\dots. $$
\end{diracproposition}

This proposition is based on Hitchin's transformation formula for spinors \cite{hitchin:74} 
(see section~\ref{compspecsec} for a proof). 

In section~\ref{systolenormsection} we will define a norm on $H^1(\torus,\mZ_2)$, the $L^2$-norm. 
This norm allows us to derive explicit lower bounds for the first eigenvalue of $D^2$ on $\torus$.
This lower bound is non-trivial if the spin structure is non-trivial.
The cosystole $\si_1(\torus,\gbel)$ can also be expressed in terms of the $L^2$-norm
  $$\si_1(\torus,\gbel):=\inf\left\{\lnorm{\al}2 \,\big|\, \al\in H^1(\torus,\mZ_2),\quad \al\neq 0 \right\}.$$
(see section~\ref{systolenormsection}, in particular Proposition~\ref{sysprop} (a)).

\begin{theorem}\label{spicotheo}
Let $(\torus,\gbel)$ be a Riemannian 2-torus with spin homomorphism $\chi$. Assume that $\Lnorm{K_\gbel}1{\torus,\gbel}<4\pi$.
Then the first eigenvalue $\mu_1$ of $D^2$ satisfies
  $$\mu_1\abel \geq {4\pi^2\,\lnorm{\chi}2^2\over \exp\left(2\cS(\Lnorm{K_\gbel}1{\torus,\gbel},\Lnorm{K_\gbel}p{\torus,\gbel}
           \abel^{1-(1/p)},p,\si_1(\torus,\gbel)^{-2}\right)},$$
The equality is attained if and only if $\gbel$ is flat.
\end{theorem}

From this theorem we will obtain two corollaries estimating $\mu_1$ in terms of the \emph{systole} $\sys$, the 
\emph{spinning systole} $\spinsys$
and the \emph{non-spinning systole} $\nssys$.
\begin{eqnarray*}
   \sys(\torus,\gbel) & := & \inf\left\{\length(\ga)\,|\,\ga\mbox{ is a non-contractible loop.}\right\}\\
   \spinsys(\torus,\gbel,\chi) & := & \inf\left\{\length(\ga)\,|\,\ga\mbox{ is a loop with }\chi([\ga])=-1.\right\}\\
   \nssys(\torus,\gbel,\chi) & := & \inf\left\{\length(\ga)\,|\,\ga\mbox{ is a non-contractible loop with $\chi([\ga])=1$}\right.\\
&&\qquad \mbox{and $[\ga]$ is a primitive element in $H_1(\torus,\mZ)$.}\bigr\}
\end{eqnarray*}
An element $\al\in H_1(\torus,\mZ)$ is called \emph{primitive} if there are no $k\in \mN$, $k\geq 2$, 
$\be\in H_1(\torus,\mZ)$ with $\al=k\cdot \be$.

\begin{corollary}\label{maincoreins}
Let $(\torus,\gbel)$ be a Riemannian 2-torus with {\bfseries non-trivial} spin homomorphism $\chi$. Assume that $\Lnorm{K_\gbel}1{\torus,\gbel}<4\pi$.
Then the first eigenvalue $\mu_1$ of $D^2$ satisfies
  $$\mu_1\abel^2 \geq {\pi^2\,\nssys(\torus,\gbel,\chi)^2\over \exp\left(2\cS(\Lnorm{K_\gbel}1{\torus,\gbel},\Lnorm{K_\gbel}p{\torus,\gbel}
           \abel^{1-(1/p)},p,{\abel\over \sys(\torus,\gbel)^2}\right)}.$$
The equality is attained if and only if $\gbel$ is flat.
\end{corollary}

\begin{corollary}\label{spinsyslowest}
Let $(\torus,\gbel)$ be a Riemannian 2-torus with {\bfseries non-trivial} spin homomorphism $\chi$. Assume that $\Lnorm{K_\gbel}1{\torus,\gbel}<4\pi$.
Then the first eigenvalue $\mu_1$ of $D^2$ satisfies
  $$\mu_1\spinsys(\torus,\gbel,\chi)^2 \geq {\pi^2\over \exp\left(4\cS(\Lnorm{K_\gbel}1{\torus,\gbel},\Lnorm{K_\gbel}p{\torus,\gbel}
           \abel^{1-(1/p)},p,{\abel\over \sys(\torus,\gbel)^2}\right)}.$$
The equality is attained if and only if
\begin{enumerate}[(a)]
\item $\gbel$ is flat, \ie $(\torus,\gbel)$ is isometric to $\mR^2/\Ga$ for a suitable lattice $\Ga$, and
\item there are generators $\ga_1, \ga_2$ for $\Ga$ statisfying 
$\ga_1\perp \ga_2$, $\chi(\ga_1)=1$ and $\chi(\ga_2)=-1$.
\end{enumerate}
\end{corollary}

Using Proposition~\ref{sysprop} and the inequalities from section~\ref{invers}
the two corollaries immediately follow from Theorem~\ref{spicotheo}.


We now turn to the Laplace operator and to the Dirac operator associated to a trivial spin structure.
We recall a well-known proposition that is the analogue of Proposition~\ref{diracvergl} for the Laplacian on surfaces
(section~\ref{compspecsec}).

\begin{laplaceproposition}
Let $M$ be a compact 2-dimensional manifold with two conformal
metrics~$\ti g$ and $g=e^{2u}\ti g$. The eigenvalues of the Laplacian on functions
corresponding to $g$
and $\ti g$ will be denoted as $0=\la_0<\la_1\leq\la_2\dots$
and \/$0=\ti\la_0<\ti\la_1\leq\ti\la_2\dots$ respectively.

Then
  $$ \la_i\,\min_{m\in M}e^{2u(m)}
     \leq \witi\la_i \leq \la_i\,  \max_{m\in M} e^{2u(m)} \qquad \forall i=1,2,\dots. $$
\end{laplaceproposition}

Together with Proposition~\ref{laplaceflatspec} and Theorem~\ref{streckabsch} we obtain

\begin{theorem}\label{laplacegesamt}
Let $(\torus,\gbel)$ be a torus conformally equivalent
to $\mR^2/\Ga$, $\vol(\mR^2/\Ga)=1$. Let 
$\Ga^*$ be the lattice dual to $\Ga$.
Let $0\leq \el_0 \leq \el_1 \leq \el_2\leq\dots$ be the sequence of lengths of $\Ga^*$, and let 
$(\la_i\,|\,i=0,1,2,\dots)$ be the spectrum of the Laplacian on functions on $(\torus,\gbel)$, then 
  $$e^{-2\osc u}\, 4\pi^2\, \el_{i}^2 \leq \la_i \, \abel \leq e^{2\osc u}\, 4\pi^2\, \el_{i}^2.$$
If $\Lnorm{K_\gbel}1{\torus,\gbel}<4\pi$, then 
  $$\osc u \leq \cS\left(\Lnorm{K_\gbel}1{\torus,\gbel},\Lnorm{K_\gbel}p{\torus,\gbel}
           \abel^{1-(1/p)},p,\si_1(\torus,\gbel)^{-2}\right).$$
\end{theorem}

Note that this theorem also provides bounds for the Laplacian on forms: 
By Poincar\'e duality the spectrum on 2-forms is the same as the spectrum on functions,
and the Laplacian on 1-forms also has the same non-zero eigenvalues, 
but each with multiplicity two.
 
The theorem implies, in particular, a lower bound on the first positive eigenvalue.

\begin{theorem}\label{lapgfl}
Let $(\torus,\gbel)$ be a Riemannian 2-torus. Assume that $\Lnorm{K_\gbel}1{\torus,\gbel}<4\pi$.
Then the first positive eigenvalue $\la_1$ of the Laplacian on functions satisfies
  $$\la_1\abel \geq {4\,\pi^2\,\si_1(\torus,\gbel)^2\over \exp\left(2\cS(\Lnorm{K_\gbel}1{\torus,\gbel},\Lnorm{K_\gbel}p{\torus,\gbel}
           \abel^{1-(1/p)},p,\si_1(\torus,\gbel)^{-2}\right)}.$$
The equality is attained if and only if $\gbel$ is flat.
\end{theorem}

\begin{corollary}\label{lapgbel}
Let $(\torus,\gbel)$ be a Riemannian 2-torus. Assume that $\Lnorm{K_\gbel}1{\torus,\gbel}<4\pi$.
Then the first positive eigenvalue $\la_1$ of the Laplacian on functions satisfies
  $$\la_1\abel^2 \geq {4\,\pi^2\,\sys(\torus,\gbel)^2\over \exp\left(2\cS(\Lnorm{K_\gbel}1{\torus,\gbel},\Lnorm{K_\gbel}p{\torus,\gbel}
           \abel^{1-(1/p)},p,{\abel\over \sys(\torus,\gbel)^2}\right)}.$$
The equality is attained if and only if $\gbel$ is flat.
\end{corollary}

\begin{remark}
Theorem~\ref{lapgfl} and Corollary~\ref{lapgbel} also hold for the 
first positve eigenvalue of $D^2$, 
if the spin structure is trivial. 
Theorem~\ref{laplacegesamt} holds for the spectrum of $D^2$, if we double the multiplicities.
\end{remark}

The structure of the paper is as follows:
In the following sections (sections~\ref{overview}--\ref{lastproofsection}) we will prove our main results. 
In section~\ref{willmoreappl}, we will apply the inequalities in Proposition~\ref{streckabsch} in order 
to obtain a lower bound on the Willmore functional.
Finally, in section~\ref{beldim} we assume that a manifold of abitrary dimension $n\geq 2$ carries two spin structures. 
We derive an upper bound for the spectra of the corresponding Dirac operators.

\section{Overview}\label{overview}

We want to obtain upper and lower bounds for the eigenvalues of
the Dirac operator and the Laplace operator on a Riemannian
2-torus $(\torus,\gbel)$.

The Clifford action of the volume element on spinors anticommutes with the Dirac operator $D$.
Thus, the spectrum of~$D$ is symmetric and is uniquely determined
by the spectrum of its square~$D^2$.
Therefore we will study the spectrum of~$D^2$
instead of the spectrum of~$D$. In the literature $D^2$ is often called the Dirac Laplacian. 

In order to prove bounds on eigenvalues we use the uniformization theorem which
tells us that we can write $\gbel$ as $\gbel=e^{2u}\gflach$ with
a real-valued function $u$ and a flat metric $\gflach$.
For flat tori the spectrum of the Laplacian and the Dirac operator is known: 
the spectra can be calculated in terms of
the dual lattice corresponding to $(\torus,\gflach)$.

We obtain bounds through the following steps.
\begin{enumerate}[(a)]
\item Comparison of the spectrum of $(\torus,\gbel)$ and the spectrum of
$(\torus,\gflach)$ (Pro\-po\-si\-tions~\ref{laplacevergl} and \ref{diracvergl}).
\item Introduction of certain spin-conformal invariants that contain information about the
dual lattice corresponding to $(\torus,\gflach)$ (section~\ref{systolenormsection}). 
\item The knowledge of spectra of flat tori (section~\ref{speczwei}).
\item A bound on $\osc u=\max u -\min u$ (section~\ref{streckabschsec}).
\item Derivation, in section~\ref{invers}, of certain inequalities that are in a sense 
inverse to the inequalities in Proposition~\ref{sysprop} and contain a curvature term.
\end{enumerate}
In section~\ref{mainproofs}, we combine the inequalities and derive the main results.

\section{Spin structures}\label{spinsection}

The eigenvalues of $D$ depend on the spin structures and 
we want to find estimates depending on the spin structure.
In this section we recall some important facts about spin structures.
Good references about spin structures are \cite{lawson.michelsohn:89}, 
\cite{bourguignon.gauduchon:92} and \cite[section II]{swift:93}. 
We will define spin structures without fixing a Riemannian metric. This definition will allow us to identify spin structures on 
diffeomorphic but not isometric manifolds (see Proposition~\ref{diracvergl}). 

Let $M$ be an oriented manifold of dimension $n\geq 2$. 
The bundle $\GLp(M)$ of oriented bases over $M$ is a principal $\GLp(n,\mR)$-bundle.
The fundamental group of $\GLp(n,\mR)$ is $\mZ$ for $n=2$ and $\mZ_2$ for $n\geq 3$. Therefore
$\GLp(n,\mR)$ has a unique connected double covering $\Th:\GLpt(n,\mR)\to\GLp(n,\mR)$.

\begin{definition}
A spin structure on $M$ is a pair $(\GLpt(M),\th)$ where $\GLpt(M)$ is a principal
$\GLpt(n,\mR)$-bundle over $M$ and $\th$ is a double covering $\GLpt(M)\to\GLp(M)$ such that
\begin{equation}\label{spincompat}
\begin{array}{cccl}
\GLpt(M) \times \GLpt(n,\mR) &\rightarrow& \GLpt(M) & \\
& & &\searrow \\
\downarrow\th\times\Theta & &   \downarrow\th & \quad M\\
& & &\nearrow \\
\GLp(M) \times \GLp(n,\mR) &\rightarrow& \GLp(M) &
\end{array}
\end{equation}
commutes. The horizontal arrows are given by the group action. 
\end{definition}

There is a spin structure on $M$ if and only if 
the second Stiefel-Whitney class $w_2(TM)$ vanishes. Such manifolds are called \emph{spin}.
From now on we assume that $M$ is spin. 

Two spin structures $(\GLpt(M),\th)$ and $(\GLpte(M),\th_1)$ are identified if there is
a fiber preserving isomorphism of principal $\GLpt(n,\mR)$-bundles 
$\al:\GLpt(M)\to \GLpte(M)$ with $\th=\th_1\circ \al$.

The set of all spin structures $(\GLpt(M),\th)$ over $M$ will be denoted by $\spinset(M)$.
The set $\spinset(M)$ has the structure of an affine space associated to the 
vector space $H^1(M,\mZ_2)$,
\ie $H^1(M,\mZ_2)$ acts freely and transitively on $\spinset(M)$.
We will describe this action: Elements in $H^1(M,\mZ_2)$ can be viewed as 
principal $\mZ_2$-bundles over $M$ \cite[Appendix~A]{lawson.michelsohn:89}. 
Let $\pi:P_\chi\to M$ be the $\mZ_2$-bundle defined by $\chi\in H^1(M,\mZ_2)$. 
Let $(\GLpt(M),\th)$ be a spin structure.
The group $\mZ_2$ acts by deck-transformation both on $\GLpt(M)$ and $P_\chi$.
We define 
  $$\GLpte(M):=(\GLpt(M)\times_M P_\chi)/\mZ_2$$
where $\mZ_2$ acts diagonally on the fiberwise product of the bundles. 
The map 
  $$\th\times_M \pi:\GLpt(M)\times_M P_\chi \to  \GLp(M)\qquad (A,\al)\mapsto \th(A)$$
is
invariant under the $\mZ_2$-action and therefore defines a map 
$\th_1:\GLpte(M)\to \GLp(M)$ compatible with \eref{spincompat}. The action of $\chi$ maps $(\GLpt(M),\th)$ to the spin structure 
$(\GLpte(M),\th_1)$. This action is free and transitive \cite[II\S 1]{lawson.michelsohn:89}.

\komment{
Note that our definition of spin structures differs from the definition given in most textbooks:
we defined spin structures without referring to a fixed Riemannian metric 
\cite[Remark~1.9]{lawson.michelsohn:89}.
Therefore we can identify spin structures on diffeomorphic manifolds that are not isometric.
This will be important in Proposition~\ref{diracvergl}. 
}

Now we fix a Riemannian metric $g$ on $M$. This reduces our structure group from $\GLp(n,\mR)$
to $\SO(n)$. The bundle of positively oriented orthonormal bases $\SO(M,g)$ is a 
principal $SO(n)$-bundle. The \emph{spin group} is defined by $\Spin(n):=\Th^{-1}(\SO(n))$ and is the unique connected 
double covering of $\SO(n)$.
A \emph{metric spin structure} is a pair $(\Spin(M,g),\th)$ where $\Spin(M,g)$ is a principal
$\Spin(n)$-bundle over $M$ and $\th$ is a double covering $\Spin(M,g)\to\SO(M,g)$
satisfying a compatibility condition analogous to \eref{spincompat}. 
For any spin structure $(\GLpt(M),\th)$ we obtain a metric spin structure $(\Spin(M,g),\th')$
by restriction:
  $$\Spin(M,g):=\th^{-1}(\SO(M,g))\qquad \th':=\res{\th}{\Spin(M,g)}.$$
Via this restriction map, the set of metric spin structures is in a natural one-to-one correspondance to $\spinset(M)$ \cite{swift:93}.

Metric spin structures are used to define spinors and the Dirac operator. 
Let $\ga_n:\Spin(n)\to\SU(\Si_n)$ be the complex spinor representation of $\Spin(n)$. This is a complex representation of dimension
$2^{[n/2]}$. It is irreducible for $n$ odd. For $n$ even, it consists of two irreducible components 
$\ga_n^+$ and $\ga_n^-$, $\ga_n^\pm: \Spin(n) \to \SU(\Si_n^\pm)$. 
The representation $\ga_n$ is not a pullback from a representation of $\SO(n)$.
The associated vector bundle $\Si M:=\Spin(M)\times_{\ga_n} \Si_n$ is called \emph{spinor bundle}
and its sections are \emph{spinors}.
The Dirac operator (see \cite{lawson.michelsohn:89} for a definition) is an elliptic operator acting on the space of smooth spinors. 

Large parts of this article will deal with the case $M=\torus$. In this case many of our definitions simplify.
Let $f:\mR^2\to \torus$ be a smooth covering map with deck transformation group $\mZ^2$ acting by translation. Then
\begin{eqnarray*}
  \tau_f:\torus \times \GLp(2) & \to & \GLp(\torus)\\ 
   (f(p),A) & \mapsto & (\pa_xf(p),\pa_yf(p))\cdot A
\end{eqnarray*}
yields a trivialization of $\GLp(\torus)$. 
\begin{definition}
The \emph{trivial spin structure on $\torus$} (with respect to $f$) is the one given by $\si_f:=(\GLpt(\torus),\th)$ with 
  $$\GLpt(\torus):=\torus \times \GLpt(2)\qquad \th:=\tau_f \circ (\id\times \Theta).$$
\end{definition}
Consider the bijection
  $$\iota_f:H^1(\torus,\mZ_2)\to \spinset(\torus),\quad \chi\mapsto \chi + \si_f.$$ 

The following proposition shows that $\iota_f$ does not depend on the choice of $f$. This will allow us 
to identify $H^1(\torus,\mZ_2)$ and $\spinset(\torus)$ via $\iota_f$.

\begin{proposition}
Let $(\GLpt(\torus),\th)$ be a spin structure on $\torus$. 
Let $\chi_f$ be the element in $H^1(\torus,\mZ_2)=\Hom_\mZ(H_1(M,\mZ),\mZ_2)$ with $\iota_f(\chi_f)=(\GLpt(\torus),\th)$.
Fix a complex structure $J$ on $T\torus$.

Then for any non-contractible smooth embedding 
$c:S^1\to\torus$ the following conditions are equivalent
\begin{enumerate}[(1)]
\item $\chi_f([c])=1$.
\item $(\dot c, J(\dot c)):S^1\to\GLp(\torus)$ lifts to  $\GLpt(\torus)$ via $\th$.
\end{enumerate} 
\end{proposition}

Characterization (2) is independent from the choice of $f$, characterization (1) is independent from the choice of $J$.
Therefore $\iota_f$ depends neither on $f$ nor $J$. 
The above proposition is an immediate consequence of the following lemma.
\begin{lemma}\label{embed}
Let $c:S^1\to\mR^2/\mZ^2$ be a non-contractible smooth embedding. Choose a lift $C:\mR\to \mR^2=\mC$, \ie $C(t)/\mZ^2=c(e^{2\pi i t})$ . 
Then
\begin{enumerate}[(1)]
\item the homology class $[c]\in H_1(\mR^2/\mZ^2,\mZ)$ is primitive, \ie not a multiple of another element in $H_1(\mR^2/\mZ^2,\mZ)$.
\item the map
  $$v(c):S^1\to S^1, e^{2\pi i t}\mapsto {\dot C(t)\over |\dot C(t)|}$$
has degree $0$.
\end{enumerate}
\end{lemma}

\proof{}
The curve $c$ can be lifted to the cylinder $Z:=\mR^2/\lan [c] \ran$. The lift will be denoted by $c^Z$. 
It is a simple closed curve generating $\pi_1(Z)$.
By Jordan's theorem about simple closed curves in $\mR^2$ we know that this curve divides $Z$ into two connected components
$Z^+$ and $Z^-$. Each of the components contains one end of the cylinder.

Let us assume that $[c]$ is not primitive, \ie $[c]=k\cdot  a$ with $k\in \mN$, $k\geq 2$ and $a\in H_1(M,\mZ)$ primitive.
The action of $a$ on $Z$ maps $Z^+$ to $Z^+$ and $Z^-$ to $Z^-$. Hence the image of $c^Z$ is mapped to itself. This contradicts $k\geq 2$.
Thus we have proven (1).

Now let $c_1:S^1\to \torus$ be another embedding, homotopic to $c$. 
A suitable lift $c_1^Z$ of $c_1$ divides
$Z^+$ into a bounded and an unbounded part. The bounded part has $c^Z$ and $c_1^Z$ as boundaries and has Euler
characteristic $0$. Therefore the Gauss-Bonnet theorem for the Euclidean metric on $Z$ yields $v(c)=v(c_1)$.
Thus the lemma only has to be checked for one representative in each primitive class. As this is trivial, (2) follows.
\qed

From now on we will identify $\spinset(\torus)$ with $H^1(\torus,\mZ_2)$ and 
$\Hom_\mZ(H_1(\torus,\mZ),\mZ_2)$.
Frequently, we will use the term ``\emph{spin homomorphism}'' instead of ``spin structure'' in order to indicate
that we regard the spin structure as an element in $\Hom_\mZ(H_1(\torus,\mZ),\mZ_2)$.

From Proposition~\ref{diracvergl} below it is clear that the trivial spin structure is the only spin structure 
such that $0$ is in the spectrum of the Dirac operator $D$. Therefore our definition of ``trivial spin structure'' 
coincides with the 
definition in section~\ref{statement}.

\begin{remark}
On oriented surfaces there is an alternative approach to define spin structures. We fix a conformal structure on $M$. 
Therefore $TM$ is complex line bundle. A \emph{line bundle spin structure} is a pair $(\Si^+ M,\th)$ of a 
complex line bundle $\Si^+ M$ and a map $\th: \Si^+ M\to TM$ satisfying
   $$\th(z\cdot q) = z^2 \cdot \th(q),\qquad \forall q\in \Si^+ M,\quad z\in \mC.$$ 
It is not hard to show that there is a natural bijection from the set of line bundle spin structures 
to the set of spin structures. For $M=\torus$, the trivial spin structure is characterized by the fact that for any non-contractible embedding 
$S^1\to\torus$ the tangent vector field $\dot{c}:S^1\to T\torus$ lifts to $\Si^+ \torus$. 
The line bundle spin structure definition is used by \cite{kusner.schmitt:p97} for example.
The Arf invariant \cite{kusner.schmitt:p97} can also be used to distinguish the trivial spin structure from the non-trivial ones.
The Arf invariant is equal to $-1$ for the trivial spin structure, and equal to $1$ for all others.
\end{remark}

\section{Comparing spectra of conformal manifolds}\label{compspecsec}

In this section we will compare Dirac and Laplace eigenvalues on 2-tori. 
We recall a proof of a well-known proposition 
(see e.g. \cite[Proposition 3.3]{dodziuk:82} for a more general version).

\begin{proposition}\label{laplacevergl}
Let $M$ be a compact 2-dimensional manifold with two conformal
metrics $\ti g$ and $g=e^{2u}\ti g$. The eigenvalues of the Laplacian on functions
corresponding to $g$
and $\ti g$ will be denoted as \/$0=\la_0<\la_1\leq\la_2\dots$
and \/$0=\ti\la_0<\ti\la_1\leq\ti\la_2\dots$ respectively.

Then
  $$ \la_i\,\min_{m\in M}e^{2u(m)}
     \leq \witi\la_i \leq \la_i\,  \max_{m\in M} e^{2u(m)} \qquad \forall i=1,2,\dots. $$
\end{proposition}

\proof{}
Let $f_0,\dots,f_i$ be eigenfunctions of $\De_{\gbel}$ to the eigenvalues $\la_0,\dots,\la_i$.
Let $U_i$ be the subspace of $V:=C^{\infty}(\torus)$ generated by $f_0,\dots,f_i$.
\noindent
We are bounding $\witi \la_i$ by the Rayleigh quotient:
  $$\witi \la _i \leq \max_{f\in U_i\ohne \{0\}} {(\De_{\ti g} f, f )_{\ti g}\over (f,f)_{\ti g}}.$$
We obtain for the numerator and the denominator:
\begin{eqnarray*}
   (\De_{\ti g} f, f )_{\ti g}&=&\int (\De_{\ti g} f) \bar f \;\dvol_{\ti g}=\int (\De_{g} f) \bar f \;\dvol_{g}\\
    &=&(\De_{g} f, f )_{g}\leq \la_i (f,f)_g
\end{eqnarray*}

  $$(f,f)_{\ti g} = \int f \bar f  \,\dvol_{\ti g}= \int f \bar f \,e^{-2u} \,\dvol_{g} \geq e^{-2\max u} \,(f,f)_g.$$
Therefore we obtain
  $$\witi \la_i \leq \la_i \, e^{2\max u}.$$
The other inequality can be proven in a completely analogous way.
\qed

There is a similar proposition for the Dirac operator.
\begin{proposition}\label{diracvergl}
Let $M$ be a compact manifold with two conformal metrics $\ti g$ and $g=e^{2u}\ti g$.
Let $D$ and $\witi D$ be the corresponding Dirac operators with respect to a common
spin structure.
We denote the eigenvalues of $D^2$ by
$\mu_1\leq\mu_2\leq\dots$ and the ones of ${\witi D}^2$
by $\witi\mu_1\leq\witi\mu_2\leq\dots$.

Then
  $$ \mu_i\,\min_{m\in M}e^{2u(m)}
     \leq \witi \mu_i \leq \mu_i\,  \max_{m\in M} e^{2u(m)} \qquad \forall i=1,2,\dots $$
\end{proposition}

\proof{}
Let $n:=\dim M$. We have
  $$\dvol_g= e^{nu} \dvol_{\ti g}.$$
There is an isomorphism of vector bundles \cite{hitchin:74}, \cite[Satz~3.14]{baum:81} or \cite[4.3.1]{hijazi:86}
\begin{eqnarray*}
\Si M & \to & \witi \Si M\\
\Psi & \mapsto & \witi\Psi
\end{eqnarray*}
over the identity $\id:M\to M$ satisfying
  $$\witi{D} (\witi \Psi)=e^u\witi{D\Psi}$$
and
  $$|\witi \Psi|=e^{{n-1\over 2}u}|\Psi|.$$

Let $(\Psi_i\,|\, i=1,2,\dots)$
be an orthonormal basis of the sections of $\Si M$ with
$\Psi_i$ being an eigenspinor of $D^2$ to the eigenvalue $\mu_i$.
The vector space spanned by $\Psi_1,\dots,\Psi_i$ will be denoted by $U_i$.

We can bound $\witi\mu_i$ by the Rayleigh quotient
  $$\witi\mu_i\leq \max_{\witi\Psi\in U_i\ohne \{0\}}
    {(\witi D \witi\Psi, \witi D\witi\Psi)_{\ti g}\over
    (\witi\Psi,\witi\Psi)_{\ti g}  }.$$
We look at the numerator and the denominator separately:
\begin{eqnarray*}
   (\witi D \witi\Psi, \witi D\witi\Psi)_{\ti g}
   &=& \int e^{2u} \lan \witi{D\Psi},\witi{D\Psi}\ran \,\dvol_{\ti g}\\
   &=& \int e^{2u+(n-1)u} \lan D\Psi,D\Psi\ran \,\dvol_{\ti g}\\
   &=& \int e^u\lan D\Psi,D\Psi\ran \,\dvol_g\\
   &\leq& (D\Psi,D\Psi)_{g} \max_{m\in M} \,e^u\\
   &\leq& \mu_i\,(\Psi,\Psi)_{g} \max_{m\in M}\, e^u\\[1.5ex]
   (\witi\Psi,\witi\Psi)_{\ti g}
   &=& \int \lan \witi \Psi, \witi \Psi \ran \,\dvol_{\ti g}\\
   &=& \int e^{-u}\,\lan \Psi,\Psi\ran\,\dvol_g\\
   &\geq & e^{-\max u} (\Psi,\Psi)_{g}
\end{eqnarray*}
Thus
  $$\witi\mu_i\leq \mu_i\max_{m\in M} e^{2u}$$
which is one of the inequalities stated in the proposition.

The other inequality can be proven in a completely analogous way.
\qed

\section{Systoles and norms on $H^1(\torus,\mZ_2)$}\label{systolenormsection}

In this section we define norms on the space of spin structures $\spinset(M)$.
These norms are strongly related to systoles.

Recall that for any compact Riemannian manifold $(M,\gbel)$, the space $H^1(M,\mR)$ carries
a natural $L^p$-norm defined to be the quotient norm of the $L^p$-norm on 1-forms
  $$\lnorm{\al}p := \inf\left\{\lnorm{\om}p \,|\, \mbox{$\om$ closed 1-form representing $\al$}\right\}.$$
For $p=\infty$ this norm is the so-called stable norm and for $p=\dim M$ it is invariant
under conformal changes of the metric.

In our special case $M=\torus$, we know that $\Ga^*=H^1(\torus,\mZ)=\Hom_\mZ(H_1(\torus,\mZ),\mZ)$
is a lattice in $H^1(\torus,\mR)$ and
that the surjective map
\begin{eqnarray}
  P: {1\over 2}\,\Ga^* & \to & \Hom_\mZ(H_1(\torus,\mZ),\mZ_2)=H^1(\torus,\mZ_2)\nonumber\\
     \al(\platz)&  \mapsto & (-1)^{2\al(\platz)}\nonumber
\end{eqnarray}
has kernel $\Ga^*$.
\begin{definition}
The $L^p$-norm on $H^1(\torus,\mZ_2)$ is the quotient norm of the $L^p$-norm
on $\Ga^*$ with respect to the quotient map  $P$, \ie for $\eta\in \Hom_\mZ(H_1(\torus,\mZ),\mZ_2)$
  $$\lnorm{\eta}p := \inf\left\{\lnorm{\al}p \,|\, \al\in {1\over 2}\Ga^*, \quad P(\al)=\eta\right\}.$$
\end{definition}
Therefore we have norms on the space of spin structures on $\torus$.
The $L^2$-norm is of particular interest as it is invariant under conformal changes and therefore it is
a spin-conformal invariant.
In the following section it will turn out that the smallest eigenvalue of $D^2$ on a flat torus with spin structure $\chi$ is
  $${4\pi^2\lnorm{\chi}2^2\over \area}.$$

Another quantity will be used for our estimate of $\osc u$ (section~\ref{streckabschsec}):
The \emph{cosystole} $\si_1$ is defined to be
  $$\si_1(\torus,\gbel):=\inf\left\{\lnorm{\al}2 \,|\, \al\in\Ga^*\ohne \{0\}\right\}.$$
For flat tori the first positive eigenvalue of the Laplacian is
  $${4\pi^2\,\si_1^2\over \area}.$$

The aim of the rest of this section is to relate the $L^2$-norms to some systolic data.

\begin{definition}
For a Riemannian 2-torus $(\torus,\gbel)$ with spin structure $\chi$ 
we define the \emph{systole} $\sys(\torus,\gbel)\in\mR$,
the \emph{spinning systole} $\spinsys(\torus,\gbel,\chi)\in\mR\cup \{\infty\}$ and
the \emph{non-spinning systole}
$\nssys(\torus,\gbel,\chi)\in\mR$ to be
\begin{eqnarray*}
   \sys(\torus,\gbel) & := & \inf\left\{\length(\ga)\,|\,\ga\mbox{ is a non-contractible loop.}\right\}\\
   \spinsys(\torus,\gbel,\chi) & := & \inf\left\{\length(\ga)\,|\,\ga\mbox{ is a loop with }\chi([\ga])=-1.\right\}\\
   \nssys(\torus,\gbel,\chi) & := & \inf\left\{\length(\ga)\,|\,\ga\mbox{ is a non-contractible loop with $\chi([\ga])=1$}\right.\\
&&\qquad \mbox{and $[\ga]$ is a primitive element in $H_1(\torus,\mZ)$.}\bigr\}
\end{eqnarray*}
An element $\al\in H_1(\torus,\mZ)$ is called \emph{primitive} if there are no $k\in \mN$, $k\geq 2$, 
$\be\in H_1(\torus,\mZ)$ with $\al=k\cdot \be$.

\end{definition}
These quantities have the following relationships
  $$\sys(\torus,\gbel)=\min\{\spinsys(\torus,\gbel,\chi),\nssys(\torus,\gbel,\chi)\}$$
  $$\sys(\widehat{\torus},\gbel)= \min\{2\cdot\spinsys(\torus,\gbel,\chi),\nssys(\torus,\gbel,\chi)\}$$
where $\widehat{\torus}$ is the covering of $\torus$ associated to $\ker \chi\subset \pi_1(\torus)$.
This covering is 2-fold for non-trivial $\chi$, and $\widehat{\torus}=\torus$ for $\chi \equiv 1$.

\begin{proposition}\label{sysprop}
Let $\gbel$ be any Riemannian metric on $\torus$ and let $\chi$ be any spin homomorphism.
There is a flat metric $\gflach$ which is conformal to $\gbel$.
This metric $\gflach$ is unique up to a multiplicative constant.

Furthermore, the following inequalities hold:
\nummerarray{a}{\displaystyle{\sys(\torus,\gbel)^2\over\area(\torus,\gbel)}\leq {\sys(\torus,\gflach)^2\over\area(\torus,\gflach)}
        ={\si_1(\torus,\gflach)^2}={\si_1(\torus,\gbel)^2}}
\nummerarray{b}{\displaystyle{\nssys(\torus,\gbel,\chi)^2\over\area(\torus,\gbel)}\leq {\nssys(\torus,\gflach,\chi)^2\over\area(\torus,\gflach)}}
  $$={4\Lnorm{\chi}2{\torus,\gflach}^2}={4\Lnorm{\chi}2{\torus,\gbel}^2}$$
\nummerarray{c}{\displaystyle{\spinsys(\torus,\gbel,\chi)^2\over\area(\torus,\gbel)}\leq{\spinsys(\torus,\gflach,\chi)^2\over\area(\torus,\gflach)}}
\nummerarray{d}{\displaystyle{\spinsys(\torus,\gflach,\chi)^2\over\area(\torus,\gflach)}\geq {1\over 4\Lnorm{\chi}2{\torus,\gflach}^2}}

{\rm (e)} For any $\eta \in H^1(\torus,\mZ_2)$ and $1\leq p\leq q \leq \infty$
\begin{eqnarray}
  \Lnorm{\eta}p{\torus,\gflach}\area(\torus,\gflach)^{-(1/p)}&=&\Lnorm{\eta}q{\torus,\gflach}\area(\torus,\gflach)^{-(1/q)}\nonumber\\
  \Lnorm{\eta}p{\torus,\gbel}\area(\torus,\gbel)^{-(1/p)}& \leq &\Lnorm{\eta}q{\torus,\gbel}\area(\torus,\gbel)^{-(1/q)}\nonumber\cr
  \Lnorm{\eta}2{\torus,\gflach} & = & \Lnorm{\eta}2{\torus,\gbel}\nonumber
\end{eqnarray}

{\rm (f)} For any $\eta \in H^1(\torus,\mZ_2)$ and $1\leq p\leq 2\leq q \leq \infty$
\begin{eqnarray*}
  \Lnorm{\eta}p{\torus,\gbel}\area(\torus,\gbel)^{\left({1\over 2}-{1\over p}\right)}
     & \leq & \Lnorm{\eta}p{\torus,\gflach}\area(\torus,\gflach)^{\left({1\over 2}-{1\over p}\right)}\\
  \Lnorm{\eta}q{\torus,\gbel}\area(\torus,\gbel)^{\left({1\over 2}-{1\over q}\right)}
     & \geq & \Lnorm{\eta}q{\torus,\gflach}\area(\torus,\gflach)^{\left({1\over 2}-{1\over q}\right)}
\end{eqnarray*}

We have equality in the inequalities of (a)--(c) if and only if $\gbel$ is flat. 

For the characterization of the equality case in (d) we choose a lattice $\Ga$ together with an isometry $I:\mR^2/\Ga\to(\torus,\gflach)$.
Then equality in (d) is equivalent to the fact 
that there are generators $\ga_1$, $\ga_2$ for the lattice $\Ga$ satisfying $\ga_1\perp \ga_2$, $I^*(\chi)(\ga_1)=1$ and 
$I^*(\chi)(\ga_2)=-1$.
\end{proposition}

\proof{}
The existence and uniqueness of $\gflach$ follows from the uniformization theorem for 2-dimensional tori.
The equations for the flat metric $\gflach$ follow directly from elementary calculations. As already stated
previously, the $L^2$-norm is invariant under conformal changes, thus the last equations in (a), (b) and (e) hold.
The inequality in (e) follows from the H\"older inequality.

The first equation in (e) then follows from the fact, that $\eta$ is represented by a real harmonic 1-form $\om$
with $\Lnorm{\eta}1{\torus,\gflach}= \Lnorm{\om}1{\torus,\gflach}$. The pointwise norm $|\om|_{\gflach}$
is constant and therefore
  $$\Lnorm{\om}1{\torus,\gflach}=\Lnorm{\om}\infty{\torus,\gflach}\area(\torus,\gflach)
    \geq \Lnorm{\eta}\infty{\torus,\gflach}\area(\torus,\gflach).$$
The inequalities in (f) follow from (e).

The remaining inequalities in (a), (b) and (c) are direct consequences from Lemma~\ref{loewnerprop} below.

The discussion of the equality case is straightforward.
\qed

\begin{lemma}[{\cite[Prop.~3.7.2]{ammanndiss}}]\label{loewnerprop}
With the notation of the previous proposition we define for $v\in H_1(\torus,\mZ)$
  $$\cL_\gbel(v):= \min \left\{\length_{\gbel}(c) \,\Bigm|\,c:S^1\to \mT^2
    \mbox{ represents }v\right\},$$
and $\cL_{\gflach}(v)$ similarly.
Then
  $${\cL_{\gbel}(v)^2\over \area(\torus,\gbel)}\leq {\cL_{\gflach}(v)^2\over \area(\torus,\gflach)}.$$
We have equality for $v\neq 0$ if and only if $\gbel$ is flat.
\end{lemma}

\proof{of lemma}
The proof of lemma follows the pattern of the proof of Loewner's theorem in \cite[4.1]{gromov:81}.

Let $\gbel=e^{2u}\gflach$. We start with a minimizer $c$ of $\cL_{\gflach}(v)$.
There is an isometric torus action on $(\torus,\gflach)$ acting by translations. Translation by $x\in \torus$
will be denoted by $L_x$. An easy calculation shows that
  $$\int_{\torus,\gflach}dx\, \length_{\gbel}(L_x(c))=  \cL_{\gflach}(v)\int_{\torus,\gflach}dx\,e^{u(x)}  \leq \cL_{\gflach}(v)\,\area(\torus,\gflach)^{1/2}\,\area(\torus,\gbel)^{1/2}.$$
Because the left hand side is an upper bound for $\cL_{\gbel}(v)\,\area(\torus,\gflach)$ the inequality of the lemma follows.
The case of equality is then obvious. \qed

\section{Spectra of flat 2-tori}\label{speczwei}

In this section we recall the well-known formulas for the spectrum of the Laplacian and of the Dirac operator
on flat 2-tori.

Because it is clear how the eigenvalues change under rescaling
we will restrict to the case
  $$\torus={\mR^2\over \Ga_{xy}} \qquad \Ga_{xy}=\Span\left\{ \pmatrix{1\cr 0},\pmatrix{x\cr y}\right\},\quad y>0$$
where $\torus$ carries the metric $\gflach$ induced by the Euclidean metric of $\mR^2$.
The dual lattice $\Gd:=H^1(\torus,\mZ)=\Hom_\mZ(\Ga_{xy},\mZ)$ is generated by the vectors
  $$\ga_1:=\pmatrix{1\cr -x/y}\quad\mbox{und}\quad\ga_2:=\pmatrix{ 0\cr 1/y}.$$
The function
\begin{eqnarray*}
  f_\ga:\torus\to \mC\quad f_\ga(x):=\exp\Big( 2\pi i\, \lan \ga,x\ran\Big)\quad \ga\in\Gd
\end{eqnarray*}
is an eigenfunction of the Laplace operator $\De$ on complex valued functions to the eigenvalue
$4\pi^2 |\ga|^2$ where $|\platz|$ denotes the Euclidean norm on $\mR^2$.
Moreover, the family $(f_\ga|\ga\in\Gd)$ is a complete system of eigenfunctions.
Note that $\ga$ can also be viewed as a $1$-form on $\torus$ and if $\lnorm{\platz}2$ is the $L^2$-norm defined in the
previous section then
  $$\lnorm{\ga}2^2=|\ga|^2\area.$$
Therefore we obtain
\begin{proposition}\label{laplaceflatspec}
The spectrum of the Laplacian on $\torus$ is given by the family
  $$\left\{{4\pi^2\lnorm{\ga}2^2\over \area}\,\Big|\,\ga\in \Gd\right\}$$
where each eigenvalue appears with the correct multiplicity.
\end{proposition}

The first three eigenvalues can be easily expressed using the invariants of the previous section
  $$\la_0=0\qquad \la_1=\la_2={4\pi^2\,\si_1^2\over \area}.$$

The eigenfunctions and eigenvalues of the square of the Dirac operator are very similar
if the spin structure is trivial. Let $\psi_1$ and $\psi_2$ be parallel orthonormal spinors
on $\torus$, then
$(f_\ga\psi_j|j=1,2;\ga\in \Gd)$ is a complete system of eigenfunctions
to the eigenvalues $4\pi^2 |\ga|^2$. Therefore the eigenvalues $\mu_1\leq \mu_2\leq \mu_3\dots$ are the same as for the
Laplace operator, but the multiplicities are doubled. In particular
  $$\mu_1=\mu_2=0\qquad\mu_3=\mu_4=\mu_5=\mu_6={4\pi^2\,\si_1^2\over \area}.$$

Now we assume that $\torus$ carries a non-trivial spin structure.
After a rescaling of the metric and an orthonormal transformation of $\mR^2$ we can assume that
the spin structure is trivial on $\pmatrix{1\cr 0}$ and non-trivial on $\pmatrix{x\cr y}$
and that
\begin{equation}\label{modulrestr}
0\leq x \leq {1\over 2}, \qquad x^2+\left(y-{1\over 2}\right)^2\geq {1\over 4}, \qquad y>0.
\end{equation}
The set of all $(x,y)$ satisfying \eref{modulrestr} is called the \emph{spin-conformal moduli space} $\spinmod$.
The elements of $\spinmod$ correspond to equivalence classes of tori with non-trivial spin structures under the equivalence relation
of conformal diffeomorphisms preserving the spin structure.

Let $(\psi_1,\psi_2)$ be a basis of parallel sections of the spinor bundle on $\mR^2$ and assume that
they are pointwise orthogonal.
Then
  $$\Psi_{j,\ga}:=\exp\Bigl( 2\pi i\,\lan \ga,x\ran\Bigr)\,\psi_j,\quad \ga\in\Gd+{\ga_2\over 2}$$
is a spinor that is invariant under the action of $\Ga_{xy}$. Thus it 
defines an eigenspinor for ${D^2:\Si\torus\to \Si\torus}$ with eigenvalue $4\pi^2 |\ga|^2$ and the family
 $(\Psi_{j,\ga}|j=1,2;\ga\in \Gd+(\ga_2/2))$ is a complete system of eigenspinors.

We obtain a similar proposition as above.
\begin{proposition}[\cite{friedrich:84}]\label{diracflatspec}
Assume that $\torus$ carries a non-trivial spin structure as above.
Then the spectrum of the square of the Dirac operator $D^2$ on $\torus$ is given by the family
  $$\left\{{4\pi^2\lnorm{\ga}2^2\over \area}\,\Big|\,\ga\in \Gd+{\ga_2\over 2}\right\},$$
and the multiplicity of each eigenvalue in the spectrum of $D^2$ is twice the multiplicity in the family.
\end{proposition}

We want to prove that $\Gd+(\ga_2/2)$ contains no vector that is shorter than $\ga_2/2$.
For this we need a lemma.

\begin{lemma}
If linearly independent vectors $v_1,v_2\in \mR^2$ satisfy
  $$0\leq \lan v_1,v_2\ran\leq|v_1|^2\leq|v_2|^2,$$
then for any integers $a,b$ with $a\neq 0$ \textbf{and} $b\neq 0$ the following inequality holds
  $$|av_1+bv_2|\geq |v_2-v_1|.$$
If $|av_1+bv_2|=|v_2-v_1|$, then $|a|=|b|=1$.
\end{lemma}

\proof{of lemma}
Let $|av_1+bv_2|\leq |v_2-v_1|$. 
Without loss of generality we can assume that $a$ and $b$ are relatively prime.
We obtain
  $$a^2|v_1|^2 - 2\,|ab|\cdot\lan v_1,v_2\ran + b^2|v_2|^2
    \leq |v_1|^2 -2\,\lan v_1,v_2\ran + |v_2|^2$$
and therefore
\begin{eqnarray*}
  (a^2+b^2-2)\,|v_1|^2&\leq& (a^2-1)\,|v_1|^2 + (b^2-1)\,|v_2|^2\\
                      &\leq& 2(|ab|-1)\,\lan v_1,v_2\ran\leq 2(|ab|-1)\,|v_1|^2.
\end{eqnarray*}
Thus $(|a|-|b|)^2\leq 0$ holds, \ie $|a|=|b|$, and as we assumed that $a$ and $b$ are relatively prime
we obtain
$|a|=|b|=1$.
Because of $|v_1+v_2|\geq |v_2-v_1|$ the lemma holds.
\qed
\begin{corollary}\ If $(x,y)\in\spinmod$, then:
\begin{enumerate}[(a)]
\item There is no vector in $\Gd+(\ga_2/2)$ that is shorter than $\ga_2/2$.
\item The shortest vectors in $\Gd\ohne \{0\}$ have length
  $$\min\left\{ {1\over y},{\sqrt{x^2+y^2}\over y}\right\}.$$
\end{enumerate}
\end{corollary}

\proof{}
\begin{enumerate}[(a)]
\item Because of relations \eref{modulrestr} the vectors $v_1:=\ga_1/2$ and $v_2:=(\ga_1+\ga_2)/2$ satisfy the conditions of the lemma.
  Any element $\ga$ of $\Gd+(\ga_2/2)$ can be written as $av_1+bv_2$, $a,b\in\mZ\ohne \{0\}$. The lemma yields
  $$|\ga|\geq |v_2-v_1|={|\ga_2|\over 2}.$$
\item This time we set $v_1=\ga_1$ and $v_2=\ga_1+\ga_2$. As before $0\leq \lan v_1,v_2\ran \leq |v_1|^2\leq |v_2|^2$.
  Any $\ga\in\Gd\ohne \{0\}$ is either a multiple of $v_1$ or $v_2$ (then $|\ga|^2\geq |v_1|^2=|\ga_1|^2=1 +(x^2/y^2))$ or
  $$|\ga|\geq |v_2-v_1|={1\over y}.$$
\end{enumerate}
\qed

Thus the smallest eigenvalue $\mu_1$ of $D^2$ satisfies
\begin{equation}\label{laom}
  \mu_1 =\pi^2|\ga_2|^2={\pi^2\over y^2}.
\end{equation}

Using the notations of the previous section we see easily that the $L^2$-norm of the spin-structure $\chi$ satisfies
  $$\lnorm{\chi}2^2={1\over 4y}.$$
With $\area=y$ we obtain
  $$\mu_1 \,\area =4\pi^2 \lnorm{\chi}2^2.$$

Analogously, we see for the cosystole that
  $$\si_1^2=\min\left\{{1\over y},{x^2+y^2\over y}\right\}.$$

\section{Regular bipartitions of 2-tori}

\begin{definition}
A \emph{regular bipartition} of $\torus$ is a pair $(X_1,X_2)$ of disjoint open subsets $X_i\subset\torus$ such that
$\pa X_1=\pa X_2$ is a smooth 1-manifold, \ie $\pa X_1=\pa X_2$ is a disjoint union of finitely many smooth circles. 
In particular
this implies $\torus=X_1\dot\cup X_2\dot\cup \pa X_1$.
\end{definition}

\begin{proposition}\label{bipartprop}
Let $(X_1,X_2)$ be a regular bipartition of $\torus$. Then exactly one of the following conditions is satisfied
\begin{enumerate}[(i)]
\item The inclusion $X_1\to \torus$ induces the trivial map $\pi_1(X_1)\to \pi_1(\torus)$.
\item The inclusion $X_2\to \torus$ induces the trivial map $\pi_1(X_2)\to \pi_1(\torus)$.
\item The boundary $\pa X_1$ has at least two components that are non-contractible in $\torus$.
\end{enumerate}
\end{proposition}

\proof{}
Assume a regular bipartition $(X_1,X_2)$ satisfies (iii), then $\pa X_1$ contains a non-contractible loop.
By a small perturbation we can achieve that this loop lies completely in $X_1$. Therefore $\pi_1(X_1)\to \pi_1(\torus)$
is not trivial. Hence $(X_1,X_2)$ does not satisfy (i). Similarly we prove that it does not satisfy (ii).

Now assume that a regular bipartition $(X_1,X_2)$ satisfies both (i) and (ii). Van-Kampen's theorem 
implies $\pi_1(\torus)=0$.
Therefore we have shown that at most one of the three conditions is satisfied.

It remains to show that at least one condition is satisfied. For this we assume that neither (i) nor (ii) is satisfied,
\ie there are continous paths $c_i:S^1\to X_i$ that are non-contractible within $\torus$.
Obviously $\pa X_1$ is homologous to zero. We will show that at least one component of $\pa X_1$ is non-homologous to zero.
Then there has to be a second component that is non-homologous to zero, because $[\pa X_1]=0$ is the sum of the homology classes
of the components.

We argue by contradiction. Assume that each component of $\pa X_1$ is homologous to zero.
Let $\pi:\mR^2\to \torus $ be the universal covering.
Then $\pi^{-1}(\pa X_1)$ is diffeomorphic to a disjoint union of countably many $S^1$.
We write
  $$\pi^{-1}(\pa X_1)=\dot{\bigcup_{i\in \mN}} Y_i$$
with $Y_i\cong S^1$.
We choose lifts $\ti c_i:\mR\to \mR^2$ of $c_i$, \ie
$\pi\left(\ti c_i(t+z)\right)=c_i(t)$ for all $t\in [0,1]$, $z\in \mZ$ and $i=1,2$.
Then we take a path $\ti\ga:[0,1]\to \mR^2$ joining $\ti c_1(0)$ to $\ti c_2(0)$.
We can assume that $\ti\ga$ is transversal to any $Y_i$.
We define $I$ to be the set of all $i\in \mN$ such that $Y_i$ meets the trace of $\ti\ga$.
The set $I$ is finite.
Using the Theorem of Jordan and Schoenfliess about simple closed curves in $\mR^2$ we can
inductively construct a compact set $K\subset \mR^2$ with boundary $\bigcup_{i\in I}Y_i$.
The number of intersections of $\ti \ga$ with $\bigcup_{i\in I}Y_i$ is odd.
Thus, either $\ti c_1(0)$ or $\ti c_2(0)$ is in the interior of $K$.
But if $\ti c_i(0)$ is in the interior of $K$, then the whole trace  $\ti c_i(\mR)$
is contained in $K$. Furthermore,
$\ti c_i(\mR)=\pi^{-1}\left(c_i([0,1])\right)$ is closed and therefore
compact. This implies that $c_i$ is homologous to zero in contradiction to our assumption.
\qed

\section{Controling the conformal scaling function}\label{streckabschsec}

Let $\torus$ carry an arbitrary metric $\gbel$.
According to the uniformization theorem we can write $\gbel=e^{2u}\gflach$ with a real function
$u:\torus\to\mR$ and a flat metric~$\gflach$.
The function $u$ is unique up to adding a constant.

The aim of this section is to estimate the quantity $\osc u:=\max u -\min u$.
The estimate is similar to an estimate of the author in a previous publication
\cite[Theorem~3.1]{ammann:00}.
The main difference is that the previous estimate needed the assumption
  $$\Lnorm{K_\gbel}p{\torus,\gbel}\left(\area(\torus,\gbel)\right)^{1-(1/p)}<4\pi$$
which is no longer needed in the estimate presented here.

\begin{theorem}\label{streckabsch}
We assume  $$\Lnorm{K_\gbel}1{\torus,\gbel}<4\pi.$$
Then for any $p\in \;\mo]1,\infty\mc[$ we obtain a bound for the oscillation of $u$
\nummerarray{a}{\osc u \leq \cS\Bigl(\Lnorm{K_\gbel}1{\torus,\gbel},\Lnorm{K_\gbel}p{\torus,\gbel}\left(\area(\torus,\gbel)\right)^{1-(1/p)},p,\si_1(\torus,\gbel)^{-2}\Bigr),}
\nummerarray{b}{\osc u \leq \cS\Bigl(\Lnorm{K_\gbel}1{\torus,\gbel},\Lnorm{K_\gbel}p{\torus,\gbel}\left(\area(\torus,\gbel)\right)^{1-(1/p)},p,{\area(\torus,\gbel)\over \sys(\torus,\gbel)^2}\Bigr),}
\par\noindent
where we use the definition
  $$\cS(\cK_1,\cK_p,p,\cV):={p\over p-1}\,\left[ {\cK_p\over 4\pi}+{1\over 2}\bigg|\log \left(1-{\cK_1\over 4\pi}\right)\bigg|+
    {\cK_1\over 8\pi -2\cK_1}\,\log\left({2\cK_p\over \cK_1}\right)\right] + {\cK_1\cV\over 8}$$
for $\cK_1>0$ and $\cS(0,\cK_p,p,\cV):=0$.
\label{schrankedef}
\end{theorem}

The function $\cS$ is continuous in $\cK_1=0$.

\begin{corollary}
Let $\cF$ be a family of Riemannian metrics conformal to the flat metric~$\gflach$.
Assume that there are constants $\cK_1\in \mo]0,4 \pi\mc[$ and $\cK_p\in \mo]0,\infty\mc[$, $p\in \;\mo]1,\infty\mc[$
with
  $$\Lnorm{K_\gbel}1{\torus,\gbel}\leq \cK_1\mbox{ and }\Lnorm{K_\gbel}p{\torus,\gbel}\left(\area(\torus,\gbel)\right)^{1-{1\over p}}
     \leq \cK_p\mbox{ \ for any }\gbel\in \cF.$$
Then the oscillation $\osc u_\gbel$ of the scaling function corresponding to $\gbel$ is uniformly bounded on $\cF$ by
  $$\osc u_\gbel\leq  \cS\Bigl(\cK_1,\cK_p,p,\cV(\torus,\gflach)\Bigr).$$
\end{corollary}

Before proving the theorem we will present some examples showing that the theorem and the corollary no longer hold
if we drop one of the assumptions
$\Lnorm{K_\gbel}1{\torus,\gbel}\leq \cK_1<4\pi$ or
$\Lnorm{K_\gbel}p{\torus,\gbel}\left(\area(\torus,\gbel)\right)^{1-{1\over p}}\leq \cK_p$.

\begin{example}
For any $\cK_1>0$ there is a sequence $(g_i)$ of Riemannian metrics with fixed conformal type, bounded volume, constant systole,
with  $$\Lnorm{K_{g_i}}1{\torus,g_i}\leq \cK_1\mbox{ \ and \ }\osc u_{g_i}\to \infty.$$
In order to construct such a sequence we take a flat torus and replace a ball by a rotationally symmetric surface which
approximates a cone for $i\to \infty$ (see \cite{ammann:00} for details).

%
%
%
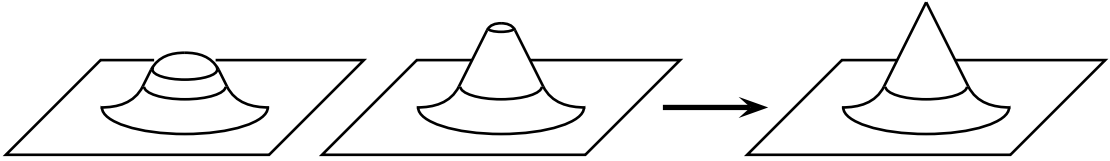
\def\flaecheneinst{\psset{linewidth=1pt,linecolor=black,linestyle=solid,fillstyle=none}}

\begin{center}
\psset{unit=.7cm}
\begin{pspicture}(-4,-1.5)(18,2)

\rput(0,0){
\flaecheneinst
\pspolygon(-3.4,-.9)(1.6,-.9)(3.4,.9)(-1.6,.9)
\psellipse(1.6,0.53)
\psframe*[linecolor=white](-1.6,0)(1.6,.53)
\psellipse(0,.4)(.8,.27)
\psframe*[linecolor=white](-.8,.4)(.8,.7)
\pspolygon*[linecolor=white](-.8,.4)(.8,.4)(0,2.0)
{\psset{unit=.8}
\psellipse(0,.9)(.8,.27)
\psframe*[linecolor=white](-.73,.9)(.73,1.2)
\psbezier(-2,0)(-1.6,0)(-1.2,0.1)(-1,0.5)
\psbezier(2,0)(1.6,0)(1.2,0.1)(1,0.5)
\psline{-}(-1,0.5)(-.8,.9)
\psline{-}(1,0.5)(.8,.9)
\psbezier(-.8,.9)(-.625,1.22)(-.32,1.3)(0,1.3)
\psbezier(.8,.9)(.625,1.22)(.32,1.3)(0,1.3)
}
}

\rput(6,0){
\flaecheneinst
\pspolygon(-3.4,-.9)(1.6,-.9)(3.4,.9)(-1.6,.9)

\psellipse(1.6,0.53)
\psframe*[linecolor=white](-1.6,0)(1.6,.53)
\psellipse(0,.4)(.8,.27)
\psframe*[linecolor=white](-.8,.4)(.8,.7)
\pspolygon*[linecolor=white](-.8,.4)(.8,.4)(0,2.0)
\psellipse(0,1.5)(.27,.088)
\psframe*[linecolor=white](-.267,1.5)(.267,1.77)
{\psset{unit=.8}
\psbezier(-2,0)(-1.6,0)(-1.2,0.1)(-1,0.5)
\psbezier(2,0)(1.6,0)(1.2,0.1)(1,0.5)
\psline{-}(-1,0.5)(-.333,1.833)
\psline{-}(1,0.5)(.333,1.833)
\psbezier(-.333,1.833)(-.266,1.966)(-.133,2)(0,2)
\psbezier(.333,1.833)(.266,1.966)(.133,2)(0,2)
}
}

\komment{
\rput(14,0){
\flaecheneinst
\pspolygon(-3.4,-.9)(1.6,-.9)(3.4,.9)(-1.6,.9)

\psellipse(1.6,0.53)
\psframe*[linecolor=white](-1.6,0)(1.6,.53)
\psellipse(0,.4)(.8,.27)
\psframe*[linecolor=white](-.8,.4)(.8,.7)
\pspolygon*[linecolor=white](-.8,.4)(.8,.4)(0,2.0)
{\psset{unit=.8}
\psellipse(0,2.1)(.2,.066)
\psframe*[linecolor=white](-.2,2.1)(.2,2.2)
\psbezier(-2,0)(-1.6,0)(-1.2,0.1)(-1,0.5)
\psbezier(2,0)(1.6,0)(1.2,0.1)(1,0.5)
\psline{-}(-1,0.5)(-.2,2.1)
\psline{-}(1,0.5)(.2,2.1)
\psbezier(-.2,2.1)(-.13,2.18)(-.08,2.2)(0,2.2)
\psbezier(.2,2.1)(.13,2.18)(.08,2.2)(0,2.2)
}}}

\psline[linewidth=2pt]{->}(9,0)(11,0)

\rput(14,0){
\flaecheneinst
\pspolygon(-3.4,-.9)(1.6,-.9)(3.4,.9)(-1.6,.9)

\psellipse(1.6,0.53)
\psframe*[linecolor=white](-1.6,0)(1.6,.53)
\psellipse(0,.4)(.8,.27)
\psframe*[linecolor=white](-.8,.4)(.8,.7)
\pspolygon*[linecolor=white](-.8,.4)(.8,.4)(0,2.0)
{\psset{unit=.8}
\psbezier(-2,0)(-1.6,0)(-1.2,0.1)(-1,0.5)
\psbezier(2,0)(1.6,0)(1.2,0.1)(1,0.5)
\psline{-}(-1,0.5)(0,2.5)
\psline{-}(1,0.5)(0,2.5)
}}
\end{pspicture}
\end{center}
\end{example}

\begin{example}
For any $\ep>0$ there is a sequence $(g_i)$ of Riemannian metrics with fixed conformal type, bounded volume, constant systole,
$-1\leq K_{g_i}\leq 1$, $\Lnorm{K_{g_i}}1{\torus,g_i}\leq 4\pi+\ep$, $\Lnorm{K_{g_i}}p{\torus,g_i}\leq \const$
and $\osc u_{g_i}\to \infty.$
In order to construct such a sequence we take a ball out of a flat torus and replace it by a hyperbolic part,
a cone of small opening angle, and a cap as indicated in the following picture. While the injectivity radius
of the hyperbolic part shrinks to zero, the oscillation of $u$ tends to infinity.
%
%
%
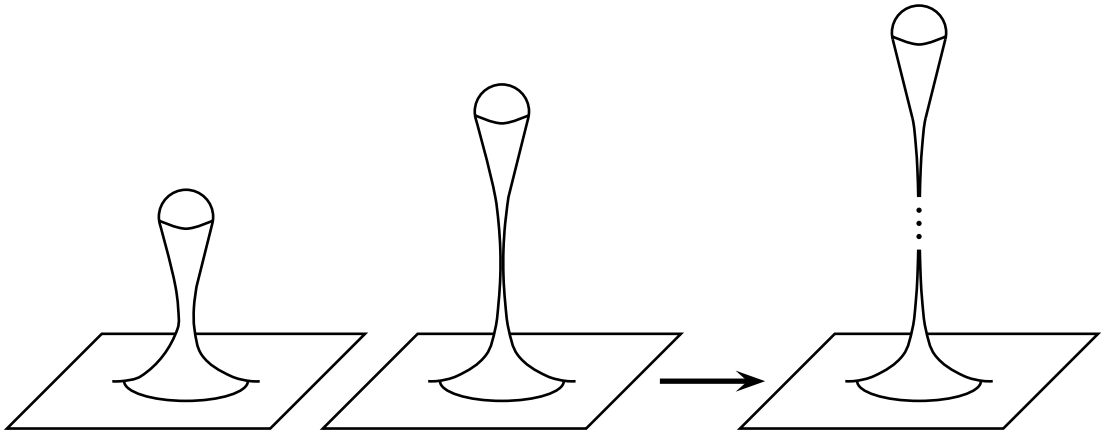
\def\flaecheneinst{\psset{linewidth=1pt,linecolor=black,linestyle=solid,fillstyle=none}}

\begin{center}
\psset{unit=.7cm}
\begin{pspicture}(-4,-1.5)(18,10)

\rput(0,0){
\flaecheneinst
\pspolygon(-3.4,-.9)(1.6,-.9)(3.4,.9)(-1.6,.9)

\psellipse(1.2,0.40)
\psframe*[linecolor=white](-1.2,0)(1.2,.40)

{\psset{unit=.5}
\rput(0,2){
\pscustom[fillstyle=solid,fillcolor=white]{
  \pscurve(2.8,-2.0)(2.5,-2.0)(1.7,-1.8)(.6,-1)(.32, 0)(.30 ,.8)(.40,1.6)
  \psarc(0,4.25){1.030776}{-14.04}{194.04}
  \pscurve(-.40,1.6)(-.30,.8)(-.32,0)(-1.7,-1.8)(-2.5,-2.0)(-2.8,-2.0)
  }
\pscurve(-1,4.1)(0,3.8)(1,4.1)
}
}
}

\rput(6,0){
\flaecheneinst
\pspolygon(-3.4,-.9)(1.6,-.9)(3.4,.9)(-1.6,.9)

\psellipse(1.2,0.40)
\psframe*[linecolor=white](-1.2,0)(1.2,.40)

{\psset{unit=.5}
\rput(0,6){
\pscustom[fillstyle=solid,fillcolor=white]{
  \pscurve(2.8,-6.0)(2.5,-6.0)(1.7,-5.8)(.6,-5)(.25,-4)(.125,-3)(.06,-2)(.06,-1)(.125,0)(.25,1)
  \psarc(0,4.25){1.030776}{-14.04}{194.04}
  \pscurve(-.25,1)(-.125,0)(-.06,-1)(-.06,-2)(-.125,-3)(-.25,-4)(-.6,-5)(-1.7,-5.8)(-2.5,-6.0)(-2.8,-6.0)
  }
\pscurve(-1,4.1)(0,3.8)(1,4.1)
}
}
}

\psline[linewidth=2pt]{->}(9,0)(11,0)

\rput(14,0){
\flaecheneinst
\pspolygon(-3.4,-.9)(1.6,-.9)(3.4,.9)(-1.6,.9)

\psellipse(1.2,0.40)
\psframe*[linecolor=white](-1.2,0)(1.2,.40)

{\psset{unit=.5}
\rput(0,9){
\psarc(0,4.25){1.03}{-14.04}{194.04}
\psline(.25,1)(1,4)
\psecurve(.6,2)(.25,1)(.125,0)(.06,-1)(.03,-2)(.015,-3)
\psline(-.25,1)(-1,4)
\psecurve(-.6,2)(-.25,1)(-.125,0)(-.06,-1)(-.03,-2)(-.015,-3)
\pscircle*(0,-2.5){.1}
\pscircle*(0,-3){.1}
\pscircle*(0,-3.5){.1}
\pscurve(-1,4.1)(0,3.8)(1,4.1)
}
\rput(0,3){
 \pscustom{
   \pscurve(2.8,-3.0)(2.5,-3.0)(1.7,-2.8)(.6,-2)(.25,-1)(.125,0)(.06,1)(.03,2)
   \gsave
      \pscurve[liftpen=1](-.03,2)(-.06,1)(-.125,0)(-.25,-1)(-.6,-2)(-1.7,-2.8)(-2.5,-3.0)(-2.8,-3.0)
      \fill[fillstyle=solid,fillcolor=white]
   \grestore}
 \pscurve(-.03,2)(-.06,1)(-.125,0)(-.25,-1)(-.6,-2)(-1.7,-2.8)(-2.5,-3.0)(-2.8,-3.0)
}
}
}
\end{pspicture}
\end{center}

In the picture the dots in the ``limit space'' indicate the hyperbolic part with 
injectivity radius tending to $0$ and diameter tending to $\infty$.
\end{example}

\proof{of Theorem~\ref{streckabsch}}
As Morse functions form a dense subset of the space of $C^\infty$-functions
with respect to the $C^\infty$-topology, we can assume without loss of generality that $u$ is a
Morse function. We set $\abel:=\area(\torus,\gbel)$ and $\afl:=\area(\torus,\gflach)$.
We define
  $$ G_<(v):=\left\{x\in \torus\,|\,u(x)<v\right\} \qquad G_>(v):=\left\{x\in \torus\,|\,u(x)>v\right\}$$
\begin{eqnarray}
  \ph:[0,\abel]& \to& \mR\nonumber\\
  A & \mapsto  & \inf \left\{\sup_{x\in X} u(x) \,\Big|\,
                X\subset \torus \mbox{ open},\; \area(X)\geq A\right\}\label{infdef}\\
      & = & \sup \left\{\inf_{x\in X^c} u(x) \,\Big|\,
                X^c\subset \torus \mbox{ open},\; \area(X^c)\geq \abel-A\right\}\label{supdef}
\end{eqnarray}

\begin{figure}[hb]
%
%
%
\newdimen\axdim
\axdim=1pt
\newdimen\cudim
\cudim=2pt
\def\ticklen{.2}
\def\abst{.5}

\begin{center}
\psset{unit=1cm}

\begin{pspicture}(-1,-1)(12,8)
\psset{linewidth=\axdim}
\psaxes[linewidth=\axdim,labels=none,ticks=none]{->}(0,0)(11.5,7)
\rput(11.9,0){$A$} 
\rput(0,7.4){$\ph(A)$}
\psline(2.4,-\ticklen)(2.4,\ticklen)
\psline(6.4,-\ticklen)(6.4,\ticklen)
\psline(8,-\ticklen)(8,\ticklen)
\psline(10,-\ticklen)(10,\ticklen)
\rput[t](2.4,-\abst){$A_-$}
\rput[t](6.4,-\abst){$A_+$}
\rput[t](8,-\abst){$A_\#$}
\rput[t](10.3,-\abst){$\matrix{\displaystyle \abel:=\hfill \cr\area(\torus,\gbel)\hfill}$}
\psline(-\ticklen,1)(\ticklen,1)
\psline(-\ticklen,2.6)(\ticklen,2.6)
\psline(-\ticklen,4.3)(\ticklen,4.3)
\psline(-\ticklen,5.7)(\ticklen,5.7)
\rput[r](-\abst,1){$\min u$}
\rput[r](-\abst,2.6){$v_-$}
\rput[r](-\abst,4.3){$v_+$}
\rput[r](-\abst,5.7){$\max u$}
\psline[linestyle=dotted]{-}(2.4,0)(2.4,2.6)
\psline[linestyle=dotted]{-}(6.4,0)(6.4,4.3)
\psline[linestyle=dotted]{-}(10,0)(10,5.7)
\psline[linestyle=dotted]{-}(0,2.6)(2.4,2.6)
\psline[linestyle=dotted]{-}(0,4.3)(6.4,4.3)
\psline[linestyle=dotted]{-}(0,5.7)(10,5.7)
\psset{linewidth=\cudim}
\psecurve[showpoints=false](-.22,1)(0,1)(.2,1.02)(1,1.8)(2,2.56)(2.4,2.6)(2.8,2.64)(3.5,2.84)(4,3.2)(5.1,3.98)(6,4.27)(6.4,4.3)(6.8,4.34)(7.4,4.6)(8,5)(8.5,5.4)(9.2,5.63)(9.8,5.69)(10,5.7)(10.2,5.7)
\end{pspicture}
\end{center}
\end{figure}

The infimum in \eref{infdef} is actually a minimum and as $u$ is a Morse function the only
minimum is attained for $X=G_<(\ph(A))$. Similarly the supremum in \eref{supdef} is attained exactly in $X^c=G_>(\ph(A))$.
The function $\ph$ is strictly increasing and is continously differentiable. The inverse of $\ph$ is given by
  $$\ph^{-1}(v)=\area(G_<(v)).$$
The differential $\ph'(A)$ is zero if and only if $\ph(A)$ is a critical value of $u$.

Now let $v\in[\min u,\max u]$ be a regular value of $u$. We obtain
\begin{eqnarray}
   \left(\ph^{-1}\right)'(v) & = & \int_{\pa G_<(v),\gbel} {1\over |du|_\gbel}
   \geq {{\length(\pa G_<(v),\gbel)}^2\over  \int_{\pa G_<(v), \gbel} |du|_\gbel}\label{phinvab}
\end{eqnarray}
where $\length(\pa G_<(v),\gbel)$ is the length of the boundary of $\pa G_<(v)$ with respect to $\gbel$.
This inequality will yield an upper bound for $\ph'$ which will provide in turn
an upper bound for $\osc u=\ph(\abel)-\ph(0)=\int_0^{\abel} \ph'$.
We transform
\begin{equation}
\int_{\pa G_<(v), \gbel} |du|_\gbel= \int_{\pa G_<(v)} *\, du =  - \int_{G_<(v),\gbel}\De_{\gbel} u
= - \int_{G_<(v),\gbel}K_\gbel.\label{gausseq}
\end{equation}
The last equation follows from the Kazdan-Warner-equation $\De_{\gbel} u=K_\gbel$ \cite{kazdan.warner:74}.
We define $\ka$ using the Gaussian curvature function
$K_\gbel:\torus\to \mR$
  $$\ka:[0,\abel]\to \mR,\quad \ka(A):=\inf \left\{\sup_{x\in X} K_\gbel(x) \,\Big|\,
    X\subset \torus \mbox{ open},\; \area(X)\geq A\right\}.$$
Any open subset $X\subset \torus$ satisfies
  $$ \int_0^{\area(X,\gbel)} \ka\leq \int_{X,\gbel} K_\gbel \leq \int_{\abel-\area(X,\gbel)}^\abel \ka$$
and for $X=\torus$ we have equality. Using Gauss-Bonnet theorem we see that
  $$\int_0^\abel \ka=0.$$
The right hand side of equation~\eref{gausseq} now can be estimated as follows.
\begin{equation}\label{kapzem}
   - \int_{G_<(\ph(A)),\gbel}K_\gbel \leq - \int_0^A \ka = \int_A^{\abel} \ka
\end{equation}
Putting \eref{phinvab}, \eref{gausseq} and \eref{kapzem} together, we obtain
  $$\ph'(A)\leq {\int_A^{\abel} \ka \over {\length(\pa G_<(\ph(A)),\gbel)}^2}.$$

Our next goal is to find suitable lower bounds for  $\length(\pa G_<(\ph(A))$.

Note that for any regular value $v$ of $u$, $(G_<(v),G_>(v))$ is a regular bipartition of $\torus$.
According to Proposition~\ref{bipartprop} exactly one of the following conditions is satisfied
\begin{enumerate}[(i)]
\item The inclusion $G_<(v)\to\torus$ induces the trivial map $\pi_1(G_<(v))\to\pi_1(\torus)$.
\item The inclusion $G_>(v)\to\torus$ induces the trivial map $\pi_1(G_>(v))\to\pi_1(\torus)$.
\item The boundary $\pa G_<(v)$ has at least two components that are non-contractible in $\torus$.
\end{enumerate}
If condition (i) is satisfied by $v$, it is obvious that it is also satisfied by $v'\in[0,v]$.
Similarly, if condition (ii) is satisfied by $v$, then it is also satisfied by $v'\in[v,\abel]$.
\begin{eqnarray*}
   v_-    & := & \sup \{v\in[0,\abel]\,|\, \mbox{(i) is satisfied for }v\}\\
   v_+    & := & \inf \{v\in[0,\abel]\,|\, \mbox{(ii) is satisfied for }v\}\\
   A_\pm  & := & \ph^{-1}(v_\pm).
\end{eqnarray*}

In each of the three cases we derive a different estimate for $\length(\pa G_<(v),\gbel)$
and therefore we obtain a different bound for $\ph'$.

\begin{enumerate}[(i)]
\item In this case $G_<(v)$ can be lifted to the universal covering $\mR^2$ of $\torus$.
We will also write $\gbel$ and $\gflach$ for the pullbacks of $\gbel$ and $\gflach$ to $\mR^2$.
The isoperimetric inequality of the flat space $(\mR^2,\gflach)$ yields
  $$\length(\pa G_<(v),\gflach)^2\geq 4\pi\, \area(G_<(v),\gflach).$$
Using the relations
\begin{eqnarray}
  \length(\pa G_<(v),\gbel)& = & e^v\, \length(\pa G_<(v),\gflach)\label{releins}\\
  \area(G_<(v),\gbel) & \leq & e^{2v}\,\area(G_<(v),\gflach)\label{relzwei}
\end{eqnarray}
we obtain
\begin{equation}
  \length(\pa G_<(v),\gbel)^2\geq 4\pi\, \area(G_<(v),\gbel).\label{negioper}
\end{equation}
Together with the H\"older inequality
  $$-\int_0^A \ka \leq\Lnorm{K_{\gbel}}{p}{\torus,\gbel} A^{1-(1/p)}$$
we get
\begin{eqnarray*}
\ph'(A) = {1 \over (\ph^{-1})'(\ph(A))} & \leq & {-\int_0^{A} \ka \over\length(\pa G_<(\ph(A)),\gbel)^2}\\
     & \leq & {1\over 4\pi}\,\Lnorm{K_{\gbel}}{p}{\torus,\gbel} A^{-{1\over p}}
\end{eqnarray*}
Integration yields
\begin{eqnarray}
  v_--\umin & = & \phi(\ph^{-1}(v_-))-\ph(0)\nonumber\\
            & \leq &  {p\over p-1}\,{1\over 4\pi}\,\Lnorm{K_{\gbel}}{p}{\torus,\gbel} (\ph^{-1}(v_-))^{1-(1/p)}\nonumber\\
            & \leq &  {p\over p-1}\,{1\over 4\pi}\,\Lnorm{K_\gbel}{p}{\torus,\gbel} (\abel)^{1-(1/p)}\label{nueins}
\end{eqnarray}
\item This case is similar to the previous one, but unfortunately because of opposite signs some
estimates do not work as before. For example \eref{relzwei} and \eref{negioper} are no longer true for
$G_<(v)$ replaced by $G_>(v)$. Instead we use Topping's inequality \cite{topping:98,topping:99}.

\begin{equation}\label{topping}
\left(\length(\pa G_>(v),\gbel)\right)^2\geq 4 \pi\, \hat A
        - 2 \int_0^{\hat A} (\hat A - a )\, \ka(\abel - a)\, da
\end{equation}
with $\hat A=\area(G_>(v),\gbel)$.
Using the estimate
\begin{eqnarray}
\int_0^{\hat A} (\hat A - a )\, \ka(\abel - a)\, da & \leq & \hat A \int _0^{\hat A} \max\{0,\ka(\abel - a)\}\, da\nonumber\\
& \leq & {\hat A \over 2}\,\Lnorm{K_\gbel}1{\torus,\gbel} \nonumber
\end{eqnarray}
we obtain
\begin{equation}
  \left(\length(\pa G_>(v),\gbel)\right)^2 \geq \left(4 \pi- \Lnorm{K_\gbel}1{\torus,\gbel}\right)\, \hat A.
\end{equation}

The obvious inequality
  $$\int_{\abel-\hat A}^\abel \ka \leq \Lnorm{\max\{0,K_{\gbel}\}}1{\torus,\gbel}\leq (1/2)\Lnorm{K_\gbel}1{\torus,\gbel}$$
yields
\begin{eqnarray*}
  \ph'(\abel-\hat A) & \leq & {1\over \hat A}\;{\Lnorm{K_\gbel}1{\torus,\gbel}\over 8\pi - 2\Lnorm{K_\gbel}1{\torus,\gbel}}.
\end{eqnarray*}

Integration yields
   $$\ph(\abel-\hat A)-\ph(A_+)\leq \log \left({\abel -A_+\over \hat A}\right)\;{\Lnorm{K_\gbel}1{\torus,\gbel}\over 8\pi - 2\Lnorm{K_\gbel}1{\torus,\gbel}}.$$
The right hand side converges to $\infty$ for $\hat A\to 0$.
Thus we have to improve our estimates for small $\hat A$.
The integral in \eref{topping} also has the following bound.
\begin{eqnarray}
\int_0^{\hat A} (\hat A - a )\, \ka(\abel - a)\, da & \leq &
    \left(\int_0^{\hat A}(\hat{A}-a)^q\,da\right)^{1/q} \cdot  \left(\int_0^{\hat A}\Big|\ka(\abel -a)\Big|^p\,da\right)^{1/p}\nonumber\\
& = & \left({{\hat A}^{q+1}\over q+1}\right)^{1/q}\cdot \Lnorm{K_\gbel}p{\torus,\gbel}
\end{eqnarray}
where we wrote $q:=p/(p-1)$ in order to simplify the notation.

We obtain a second lower bound on the length
\begin{equation}\label{laengelowzwei}
  \left(\length(\pa G_>(v),\gbel)\right)^2\geq 4 \pi\, \hat A - c {\hat A}^{1+{1\over q}}\Lnorm{K_\gbel}p{\torus,\gbel}
\end{equation}
for any  $c\geq 2/ \sqrt[q]{q+1}$, \eg $c=2$.
Note that our assumption $\Lnorm{K_\gbel}1{\torus,\gbel}<4\pi$ does not imply that the right hand side of the above
inequality is always positive.
Although \eref{laengelowzwei} is better for small $\hat A$, it is not strong enough to control the length for larger $\hat A$.
However, for
  $$\hat A < \left({4\pi\over c\cdot\Lnorm{K_\gbel}p{\torus,\gbel}}\right)^{q}$$
we use \eref{laengelowzwei} and
  $$\int_{\abel-\hat A}^\abel \ka \leq {\hat A}^{1/q}\, \Lnorm{K_{\gbel}}p{\torus,\gbel}$$
to obtain the estimate
\begin{eqnarray*}
  \ph'(\abel-\hat A) & \leq &{{\hat A}^{-1/p}\,\Lnorm{K_\gbel}p{\torus,\gbel}\over 4\pi -
  c {\hat A}^{1/q}\Lnorm{K_\gbel}p{\torus,\gbel}}.
\end{eqnarray*}
\def\lokpnorm{\Lnorm{K_\gbel}p{\torus,\gbel}}
With the substitution
  $$w=w(A)=4\pi - c(\abel -A)^{1/q}\lokpnorm$$
integration yields
\begin{eqnarray}
  \ph(\abel) - \ph(A_\#) & = & \int_{A_\#}^\abel \ph'(A)\,dA \nonumber\\
               & \leq & \int_{w(A_\#)}^{w(\abel)} {q\over c}\,{1\over w}\,dw =  {q\over c}\, \log {w(\abel)\over w(A_\#)}\nonumber\\
               & = & {q\over c} \log { 4\pi \over 4\pi -c (\abel-A_\#)^{1/q} \lokpnorm}\nonumber
\end{eqnarray}
for any $A_\#$ between $\abel-\left(4\pi/(c\cdot\Lnorm{K_\gbel}p{\torus,\gbel})\right)^{q}$ and $\abel$.
We choose
  $$A_\#:=\max\left\{\abel - \left({\Lnorm{K_\gbel}1{\torus,\gbel}\over 2\, \lokpnorm}\right)^q,A_+\right\}.$$
Finally we obtain the estimates
\begin{eqnarray}
  \umax - \ph(A_\#) & \leq & {q\over c}\, \log{8\pi \over 8\pi -c \Lnorm{K_\gbel}1{\torus,\gbel}}\label{nuzweia}\\
  \ph(A_\#) - v_+  & \leq & q \,{\Lnorm{K_\gbel}1{\torus,\gbel}\over 8\pi -2 \Lnorm{K_\gbel}1{\torus,\gbel}}\,
                       \log \left({2\,\abel^{1/q} \Lnorm{K_\gbel}p{\torus,\gbel}\over \Lnorm{K_\gbel}1{\torus,\gbel}}\right).\label{nuzweib}
\end{eqnarray}
For $c=2$ the right hand sides of these inequalities contribute
two summands to the formula for $\cS$.
\item
If $v=\ph(A)$ is a regular value of $u$ between $v_-$ and $v_+$, then $\pa G_<(v)$ contains at least two components
that are non-contractible in $\torus$. Hence, for any metric $\ti g$ on $\torus$ we get
  $$\length(\pa G_<(v),\ti g)\geq 2 \,\sys(\torus,\ti g).$$
In order to prove (a) of Theorem~\ref{streckabsch} we apply this equation to $\ti g:=\gflach$.
Using $\int_A^\abel \ka \leq (1/2) \Lnorm{K_\gbel}1{\torus,\gbel}$ and $\length(\pa G_<(v),\gbel)=e^v\,\length(\pa G_<(v),\gflach)$
we obtain
\begin{eqnarray}
  \ph'(A)& \leq & e^{-2\ph(A)} {\int_A^\abel \ka\over4\,\sys(\torus,\gflach)^2}\nonumber\\
         & \leq & {1\over 8}\, e^{-2\ph(A)}\, {\Lnorm{K_\gbel}1{\torus,\gbel}\over \sys(\torus,\gflach)^2}.
\end{eqnarray}
Integration yields
\begin{eqnarray}
v_+-v_- & = & \int_{A_-}^{A_+} \ph'(A)\,dA\nonumber\\
        & \leq & {1\over 8}\, {\Lnorm{K_\gbel}1{\torus,\gbel}\over \sys(\torus,\gflach)^2}\,\int_{A_-}^{A_+} e^{-2\ph(A)}\,dA\nonumber\\
        & \leq & {1\over 8}\, {\Lnorm{K_\gbel}1{\torus,\gbel}\over \sys(\torus,\gflach)^2}\,\afl
\end{eqnarray}
where we used $\afl=\area(\torus,\gflach)=\int_0^\abel e^{-2\ph(A)}\,dA$.

Together with inequalities \eref{nueins}, \eref{nuzweia} and \eref{nuzweib} we obtain statement (a) of the theorem.

Similarly, setting $\ti g:=\gbel$ we get statement (b).
\end{enumerate}
\qed

\section{Some ``inverse'' inequalities}\label{invers}

In Proposition~\ref{sysprop} and Lemma~\ref{loewnerprop} we proved
some inequalities relating the metric $\gbel$ to $\gflach$. It is easy to prove that they also hold in the
other direction if we add a factor like $e^{2\osc u}$.

Explicitely we obtain:
\nummerarray{a}{\displaystyle{\sys(\torus,\gbel)^2\over\area(\torus,\gbel)}\geq e^{-2\osc u}{\sys(\torus,\gflach)^2\over\area(\torus,\gflach)}}
\nummerarray{b}{\displaystyle{\nssys(\torus,\gbel,\chi)^2\over\area(\torus,\gbel)}\geq e^{-2\osc u}{\nssys(\torus,\gflach,\chi)^2\over\area(\torus,\gflach)}}
\nummerarray{c}{\displaystyle{\spinsys(\torus,\gbel,\chi)^2\over\area(\torus,\gbel)}\geq e^{-2\osc u}{\spinsys(\torus,\gflach,\chi)^2\over\area(\torus,\gflach)}}

{\rm (d)} For any $\eta \in H^1(\torus,\mZ_2)$ and $1\leq p\leq 2\leq q \leq \infty$
\begin{eqnarray*}
  \Lnorm{\eta}p{\torus,\gbel}\area(\torus,\gbel)^{\left({1\over 2}-{1\over p}\right)}
     & \geq & e^{\left(1-{2\over p}\right)\,\osc u}
     \Lnorm{\eta}p{\torus,\gflach}\area(\torus,\gflach)^{\left({1\over 2}-{1\over p}\right)}\\
  \Lnorm{\eta}q{\torus,\gbel}\area(\torus,\gbel)^{\left({1\over 2}-{1\over q}\right)}
     & \leq & e^{\left({1}-{2\over q}\right)\,\osc u}
     \Lnorm{\eta}q{\torus,\gflach}\area(\torus,\gflach)^{\left({1\over 2}-{1\over q}\right)}
\end{eqnarray*}

A combination of these inequalities together with our upper bound for $\osc u$ in the previous section
enables us to compare the quantities under consideration for a flat and an arbitrary metric in the same (spin-)conformal
class.

\section{Proof of the main results}\label{mainproofs}\label{lastproofsection}

Combining the inequalities derived in the previous sections, we are now able to derive our main results.

Theorem~\ref{maineins} is a consequence of Proposition~\ref{diracflatspec} together with Proposition~\ref{diracvergl} and Theorem~\ref{streckabsch}.
Theorem~\ref{spicotheo} then follows from the calculation of the first eigenvalue of $D^2$ on flat tori at the 
end of section~\ref{speczwei}. Using the inequalities in Proposition~\ref{sysprop} and section~\ref{invers} we can derive
Corollaries~\ref{maincoreins} and \ref{spinsyslowest}.  

Similarly, Theorem~\ref{laplacegesamt} is a consequence of Proposition~\ref{laplaceflatspec} together with  Proposition~\ref{laplacevergl} and 
Theorem~\ref{streckabsch}.
Theorem~\ref{lapgfl} then follows from the calculation of the first positive eigenvalue of $\De$ on flat tori at the 
end of section~\ref{speczwei}. Using the inequalities in Proposition~\ref{sysprop} we obtain
Corollary~\ref{lapgbel}.

\section{An application to the Willmore functional}\label{willmoreappl}

In this section $S^3$ always carries the metric $g_{S^3}$ of constant sectional curvature $1$.
For any immersion $F:\torus\to S^3$ we define the Willmore functional
  $$\cW(F):=\int_{(\torus,F^*g_{S^3})} |H_{\torus\to S^3}|^2+1$$
where $H$ is the relative mean curvature of $F(\torus)$ in $S^3$ and integration is the usual integration of functions $\torus\to\mR$
over the Riemannian manifold $(\torus,F^*g_{S^3})$. 
Note that the mean curvature $H$ of $F(\torus)$ in $\mR^4$ satisfies
  $$|H|^2=|H_{\torus\to S^3}|^2+1.$$

The Willmore conjecture states that $\cW(F)\geq 2\pi^2$. Li and Yau \cite[Fact~3]{li.yau:82}
proved that the conjecture holds if $F$ is not an embedding.

Any immersion $F:\torus\to S^3$ induces a spin structure $\ph_F$ on $\torus$.
The spin structure $\ph_F$ is non-trivial if and only if
$F$ is regularly homotopic to an embedding.
Thus for any immersion $F$ which is regularly homotopic to an embedding, the pair $(F^*g_{S^3},\ph_F)$ 
defines an element $(x,y)$ in the spin-conformal moduli
space $\spinmod$ (defined in section~\ref{speczwei}). 
In order to shorten our notation we write $[F]:=(x,y)\in\spinmod$.
If $\torus$ already carries a spin structure, we say that $F$ is \emph{spin} iff $\ph_F=\ph$.

Li and Yau proved:
\begin{theorem}[{\cite[Theorem~1]{li.yau:82}}]
Let $F:(\torus,g)\to (S^3,g_{S^3})$ be a conformal embedding, let $\abel$ be the area of $(\torus,g)$ and let
$\la_1$ be the first positive eigenvalue of the Laplacian $\De$ on $(\torus,\gbel)$ then
  $$\cW(F)\geq {1\over 2}\,\la_1\abel.$$
\end{theorem}

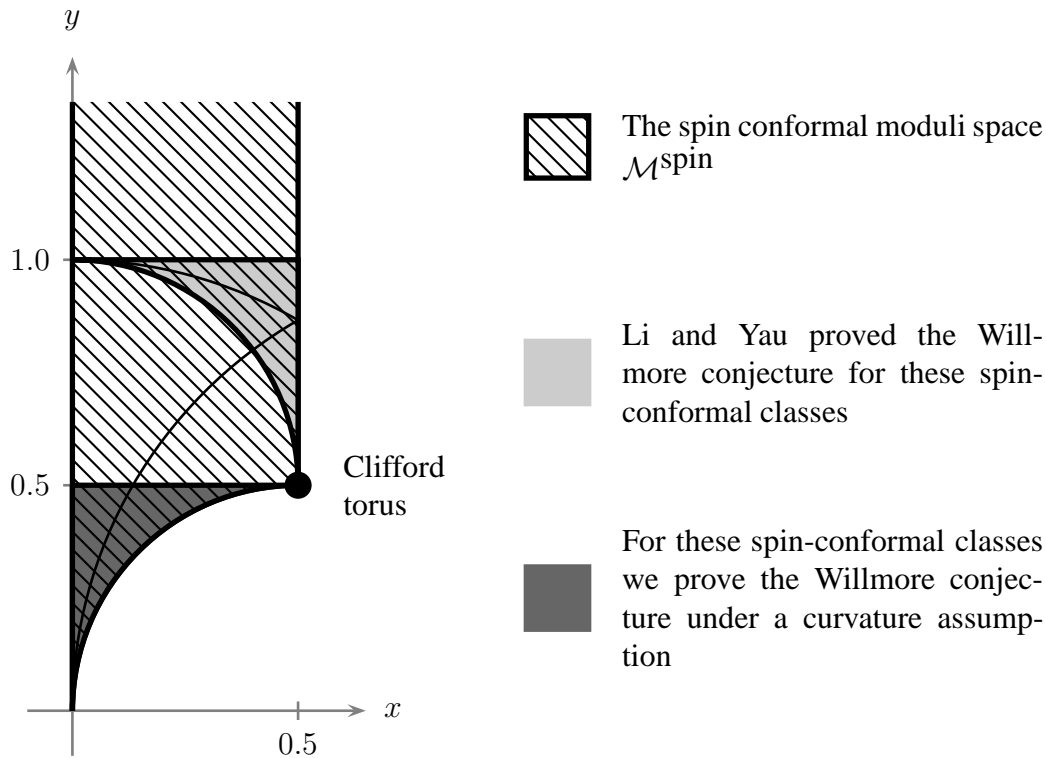
\begin{figure}[tb]
%
%
%

\newgray{abshilflinfarbe}{.30}
\newgray{absliyaufarbe}{.80}  
\newgray{absammannfarbe}{.40}  

\newdimen\descwidth
\descwidth=5.6cm
\def\achseneinst{\psset{linewidth=1pt,linecolor=gray,linestyle=solid,fillstyle=none}}

\def\abspinpfeileinst{\psset{linewidth=1pt,linecolor=black,linestyle=solid,fillstyle=none}}
\def\absmodformeinst{\psset{linewidth=2pt,linecolor=black,linestyle=solid,fillstyle=vlines,hatchcolor=black}}
\def\absliyaueinst{\psset{linewidth=2pt,linecolor=black,linestyle=solid,fillstyle=vlines*,fillcolor=absliyaufarbe,hatchcolor=black}}
\def\absammanneinst{\psset{linewidth=2pt,linecolor=black,linestyle=solid,fillstyle=vlines*,fillcolor=absammannfarbe,hatchcolor=black}}
\def\absralilineinst{\psset{linewidth=1pt,linecolor=black,linestyle=solid,fillstyle=none}}
\def\cliffordeinst{\psset{linewidth=2pt,linecolor=black,linestyle=solid,fillstyle=solid,fillcolor=black}}

\psset{unit=3cm}
\begin{pspicture}(-.2,-.2)(4.2,3.35)
\achseneinst
\psaxes[Dx=.5,dx=1,Dy=.5,dy=1]{->}(0,0)(-.2,-.2)(1.3,2.9)

\absmodformeinst
\pscustom{
  \psline(0,2.7)(0,0)
  \psarcn(1,0){1}{180}{90}
  \psline(1,1)(1,2.7)
} 

\absliyaueinst
\pscustom{
  \psarc(0,1){1}{0}{90}
  \psline(0,2)(1,2)
  \psline(1,2)(1,1)
}

\absammanneinst
\pscustom{
  \psline(0,1)(0,0)
  \psarcn(1,0){1}{180}{90}
  \psline(1,1)(0,1)
} 

\absralilineinst
\psarc(0,0){2}{60}{90}
\psarc(2,0){2}{120}{180}

\cliffordeinst
\pscircle(1,1){.06}

\rput(1.42,0){$x$}
\rput(0,3.06){$y$}

\rput[l](1.2,1){\vbox{Clifford\\torus}}

\def\tere{2.0}
\def\terz{2.4}
\rput[r](\tere,1.5){\psframe[linecolor=absliyaufarbe,fillstyle=solid,fillcolor=absliyaufarbe](0,-.15)(.3,.15)}
\rput[l](\terz,1.5){ 
\vtop{\hsize=\descwidth\noindent Li and Yau proved the Willmore conjecture for these spin-conformal classes}
}
\rput[r](\tere,.5){\psframe[linecolor=absammannfarbe,fillstyle=solid,fillcolor=absammannfarbe](0,-.15)(.3,.15)}
\rput[l](\terz,.5){
\vtop{\hsize=\descwidth\noindent For these spin-conformal classes we prove the Willmore conjecture under a curvature assumption} 
}
\absmodformeinst
\rput[r](2.0,2.5){\psframe(0,-.15)(.3,.15)}
\rput[l](2.4,2.5){
\vtop{\hsize=\descwidth\noindent 
The spin conformal moduli space $\spinmod$}
}
\end{pspicture}
\caption{The spin conformal moduli space}\label{spinmodulbild}
\end{figure}

From this theorem the conjectured inequality $\cW(F)\geq 2\pi^2$ follows, if $[F]$ lies in
a compact subset of $\spinmod$ with positive measure
(see Figure~\ref{spinmodulbild}).

A similar lower bound for $\cW(F)$ in terms of Dirac eigenvalues has been given by B\"ar.
\begin{theorem}[{\cite{baer:98}}]\label{baertorabsch}
Let $F:(\torus,g,\ph)\to (S^3,g_{S^3})$ be an {\bfseries isometric} spin immersion.
Then for the first eigenvalue $\mu_1$ of the square of the Dirac operator the inequlity
  $$\cW(F)\geq \mu_1\abel$$
holds.
\end{theorem}
Note that this estimate is only non-trivial if $F$ is regularly homotopic to an embedding.
\begin{remark}
At the end of this section we will show by example that in general ``isometric'' can not be replaced be ``conformal'' in this theorem.
\end{remark}

Our goal now is to apply our previous estimates and derive lower bounds for $\cW(F)$. One way to deduce such bounds is   
to combine the theorem with our lower estimates for the first eigenvalue of the square of the Dirac operator. 
These lower etimates for the Willmore functional are weaker than the ones derived by the author in \cite{ammann:00},
therefore we skip this approach.

In this article, our approach is to modify the techniques of Theorem~\ref{baertorabsch}. This yields
together with Theorem~\ref{streckabsch} new results
about the Willmore functional. 

As in the previous sections we define
  $$\cS(\cK_1,\cK_p,p,\cV):={p\over p-1}\,\left[ {\cK_p\over 4\pi}+{1\over 2}\bigg|\log \left(1-{\cK_1\over 4\pi}\right)\bigg|+
    {\cK_1\over 8\pi -2\cK_1}\,\log\left({2\cK_p\over \cK_1}\right)\right] + {\cK_1\cV\over 8}$$
for $\cK_1>0$ and $\cS(0,\cK_p,p,\cV):=0$.

\begin{theorem}\label{addabschsatz}
Let $F:\torus\to S^3$ be an immersion of the 2-dimensional torus in $S^3$ carrying the standard metric $g_{S^3}$.
Let $F$ be regularly homotopic to an embedding. We set $\gbel:=F^*g_{S^3}$. Let $(x,y)=[F]\in \spinmod$.
Then
  $$\will(F)\geq {\pi^2\over y}- {1\over 8}\,(\osc u)\,\left\|K_{\gbel}\right\|_{L^1(\torus,\gbel)}.$$
In particluar if $\Lnorm{K_{\gbel}}1{\torus,\gbel}<4\pi$ and any $p>1$
  $$\will(F)\geq {\pi^2\over y}- {1\over 8}\,\cS
         \left\|K_{\gbel}\right\|_{L^1(\torus,\gbel)}$$
with $\cS:=\cS\left(\Lnorm{K_\gbel}1{\torus,\gbel},\Lnorm{K_\gbel}p{\torus,\gbel}
           \abel^{1-(1/p)},p,{\abel\over \sys(\torus,\gbel)^2}\right)$ or\newline
$\cS:=\cS\left(\Lnorm{K_\gbel}1{\torus,\gbel},\Lnorm{K_\gbel}p{\torus,\gbel}
           \abel^{1-(1/p)},p,\si_1(\torus,\gbel)^{-2}\right)$.
\end{theorem}

\proof{}
We write the induced metric $\gbel$ on $\torus$
in the form $\gbel =e^{2u}\gflach$ with $\gflach$ flat.
Any Killing spinor on $S^3$ with the Killing constant $\al=(1/2)$ induces a spinor $\psi$
on $(\torus,\gbel)$ satisfying
  $$D_{\gbel}\psi= H \psi + \nu \psi,$$
where
  $$\nu=\ga(e_1)\ga(e_2)=\pmatrix {-i & 0 \cr 0& i }\in \End\left(\Si^{+}\torus \oplus \Si^{-}\torus\right)$$
(see \eg \cite{baer:98}).
There is an isomorphism of vector bundles \cite{hitchin:74},\cite[4.3.1]{hijazi:86}
\begin{eqnarray*}
  \Si \torus& \to & \witi\Si \torus\\
  \Psi &\mapsto \wihat \Psi
\end{eqnarray*}
with
  $$e^u \,\widehat{\overbrace{D_{\gbel}\Psi}} = D_{\gflach}\wihat\Psi + {1\over 2}\ga_{\gflach}(\grad_{\gflach} u) \wihat\Psi$$
and
  $$|\wihat\Psi|= |\Psi|.$$
Here $\ga_{\gflach}$ means Clifford multiplication corresponding to the metric~$\gflach$.
Note that $\witi\Psi$ from section~\ref{compspecsec} satisfies $\witi\Psi=e^{(u/2)}\wihat\Psi$.

We apply this transformation for $\Psi=\psi$ and we obtain
  $$D_{\gflach}\wihat\psi=-{1\over 2}\ga_{\gflach}(\grad_{\gflach} u) \wihat\psi + e^u H \wihat\psi
    + e^u \nu \wihat\psi.$$
As $\nu$, $\ga(V)$ and $\nu\ga(V)$ are skew-hermitian for any vector $V$, this yields
\begin{eqnarray*}
  \left|D_{\gflach}\wihat\psi\right|^2
     & = & {1\over 4}\left|\ga_{\gflach}(\grad_{\gflach} u) \wihat\psi\right|^2
           + e^{2u} H^2 \left|\wihat\psi\right|^2
           + e^{2u} \left|\nu \wihat\psi\right|^2\\
     & = & {1\over 4}\left|du\right|_{\gflach}^2 + e^{2u} H^2 + e^{2u}.
\end{eqnarray*}

Integration over $(\torus,\gflach)$ provides
  $$\ti\la_1\area(\torus,\gflach)\leq {1\over 4}\int_{\torus}\left|du\right|_{\gflach}^2\,\dvol_{\gflach}
    + \will(F),$$
where $\ti\la_1$ denotes the smallest eigenvalue of the square of the Dirac operator on $(\torus,\gflach)$.

On the other hand
\begin{eqnarray*}
   \int_{\torus}\left|du\right|_{\gflach}^2\dvol_{\gflach}
   & = & \int_{\torus} u\De_{\gflach} u\,\dvol_{\gflach} \\
   & = & \int_{\torus} e^{2u}u K_g\,\dvol_{\gflach}\\
   & = & \int_{\torus} u K_g\,\dvol_{\gbel}\\
   &\leq& {1\over 2}\,(\osc u)\,\|K_g\|_{L^1(\torus,\gbel)}.
\end{eqnarray*}
Together with Theorem~\ref{streckabsch} and the results of section~\ref{speczwei} we get the statement.
\qed

\begin{corollary}
For any $\ka_1\in\mo]0,4\pi\mc[$, any $p>1$ and any $\ka_p>0$ there is a neighborhood $U$ of the $(y\to 0)$-end of $\spinmod$
with the following property:\newline
If $F:\torus\to S^3$ is an immersion such that the induced metric $\gbel:=F^*g_{S^3}$ and the
induced spin structure $\phi_F$ represent a spin-conformal class in $U$ and if the curvature conditions
  $$\Lnorm{K_\gbel}1{\torus,\gbel}<\ka_1 \quad\mbox{and}\quad \Lnorm{K_\gbel}p{\torus,\gbel}\abel^{1-(1/p)}<\ka_p$$
are satisfied,
then the Willmore conjecture
  $$W(F)\geq 2\pi^2$$
holds.
\end{corollary}

\begin{corollary}\label{willmoreinfty}
Let $F_i:\torus\to S^3$ be a sequence of immersions. The induced metrics $g_i:=F_i^*g_{S^3}$ together with the induced spin structures define
a sequence $(x_i,y_i)$ in the spin-moduli space $\spinmod$. Assume that $y_i\to 0$ and that the curvature conditions
  $$\Lnorm{K_{g_i}}1{\torus,g_i}<\ka_1<4\pi \quad\mbox{and}\quad \Lnorm{K_\gbel}p{\torus,\gbel}\abel^{1-(1/p)}<\ka_p$$
are satisfied for some $p>1$ and  $\ka_p<\infty$.
Then
  $$\cW(F_i)\to \infty.$$
\end{corollary}

The conclusion of the second corollary is false if we drop the curvature conditions. 
To see this we construct
a sequence of immersions with $y_i\to 0$ and $\cW(F_i)<\const$.
We start with an embedding $F:\torus\to S^3$ which looks in a neighborhood of
some point like a cylinder. Now we ``strangle'' the torus as in the picture below:
%
%
%
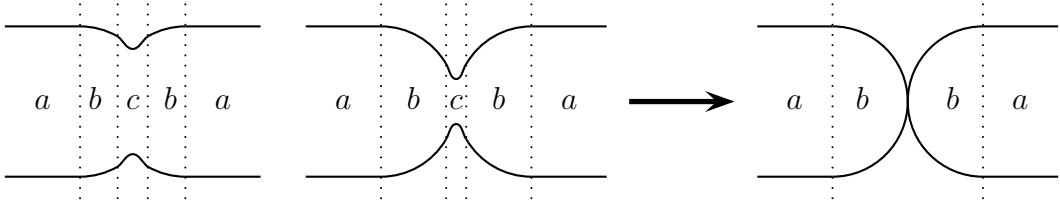
\def\flaecheneinst{\psset{linewidth=1pt,linecolor=black,linestyle=solid,fillstyle=none}}

\begin{center}
\psset{unit=1cm}
\begin{pspicture}(1,-1)(15,1)

\rput(0,0){
     \psline(1,1)(2,1)
     \psarcn(2,0){1}{90}{60}
     \psecurve(2.0,1.15)(2.5,.866)(2.7,.7)(2.9,.866)(3.4,1.15)
     \psarcn(3.4,0){1}{120}{90}
     \psline(3.4,1)(4.4,1)
     \psline(1,-1)(2,-1)
     \psarc(2,0){1}{270}{300}
     \psecurve(2.0,-1.15)(2.5,-.866)(2.7,-.7)(2.9,-.866)(3.4,-1.15)
     \psarc(3.4,0){1}{240}{270}
     \psline(3.4,-1)(4.4,-1)
     \psline[linestyle=dotted](2,-1.3)(2,1.3)
     \psline[linestyle=dotted](2.5,-1.3)(2.5,1.3)
     \psline[linestyle=dotted](2.9,-1.3)(2.9,1.3)
     \psline[linestyle=dotted](3.4,-1.3)(3.4,1.3)
     \rput[b](1.5,-.1){$a$}
     \rput[b](3.9,-.1){$a$}
     \rput[b](2.2,-.1){$b$}
     \rput[b](3.2,-.1){$b$}
     \rput[b](2.7,-.1){$c$}
}

\rput(4,0){
     \psline(1,1)(2,1)
     \psarcn(2,0){1}{90}{30}
     \psecurve(2.366,1.366)(2.866,.5)(3.0,.3)(3.133,.5)(3.633,1.15)
     \psarcn(4,0){1}{150}{90}
     \psline(4,1)(5,1)
     \psline(1,-1)(2,-1)
     \psarc(2,0){1}{270}{330}
     \psecurve(2.366,-1.366)(2.866,-.5)(3.0,-.3)(3.133,-.5)(3.633,-1.15)
     \psarc(4,0){1}{210}{270}
     \psline(4,-1)(5,-1)
     \psline[linestyle=dotted](2,-1.3)(2,1.3)
     \psline[linestyle=dotted](2.866,-1.3)(2.866,1.3)
     \psline[linestyle=dotted](3.133,-1.3)(3.133,1.3)
     \psline[linestyle=dotted](4,-1.3)(4,1.3)
     \rput[b](1.5,-.1){$a$}
     \rput[b](4.5,-.1){$a$}
     \rput[b](2.433,-.1){$b$}
     \rput[b](3.566,-.1){$b$}
     \rput[b](3,-.1){$c$}
}

\psline[linewidth=2pt]{->}(9.3,0)(10.7,0)

\rput(10,0){
     \psline(1,1)(2,1)
     \psarcn(2,0){1}{90}{270}
     \psline(2,-1)(1,-1)
     \psline(5,-1)(4,-1)
     \psarcn(4,0){1}{270}{90}
     \psline(4,1)(5,1)
     \psline[linestyle=dotted](2,-1.3)(2,1.3)
     \psline[linestyle=dotted](4,-1.3)(4,1.3)
     \rput[b](1.5,-.1){$a$}
     \rput[b](4.5,-.1){$a$}
     \rput[b](2.4,-.1){$b$}
     \rput[b](3.6,-.1){$b$}
}
\end{pspicture}
\end{center}

We get a sequence $F_i:\torus\to S^3$  of $C^1$-embeddings with the following properties:
\begin{enumerate}[(i)]
\item $F_i(\torus)$ coincides with $F(\torus)$ in region $a$  
\item $F_i(\torus)$ coincides with a part of a half-sphere in region $b$,
\item $F_i(\torus)$ coincides with a minimal surface in region  $c$.
\end{enumerate}
Note that the regions $a$, $b$ and $c$ depend on $i$. In the limit $i\to \infty$,
region $c$ disappears.
After smoothing we get a family of smooth embeddings satisfying  
both $y_i\to 0$ and $\cW(F_i)<\const$ and $\area(\torus, F_i^* g_{S^3})\to\const$.

Hence, the first eigenvalue of $D^2$ is bounded from above. But the
first eigenvalue of the spin-conformally equivalent flat torus with unit volume 
converges to $\infty$. This implies that there are spin-conformal classes
in which the optimal constants in Lott's inequality \eref{lott} 
are not attained by flat metrics.

From this example we can also conclude that Theorem~\ref{baertorabsch}
does no longer hold, if we replace the condition ``isometric spin immersion'' by 
``conformal spin immersion''.

\section{Comparing spectra for different spin structures}\label{beldim}

\def\homgitter#1{H_{\mR}^1\left(M,#1\mZ\right)}

In this section we remove the assumption $\dim M=2$ and assume that a compact Riemannian 
spin manifold~$(M,g)$ of arbitrary dimension
carries at least two different spin structures $\th$ and $\th'$.
The space of spin structures on $M$ is an affine space associated to the vector
space $H^1(M,\mZ_2)$ which will be identified with $\Hom_{\mZ}(H_1(M,\mZ),\mZ_2)$
and $\Hom(\pi_1(M),\mZ_2)$.

For $r\in \mR$, let $\homgitter{r}$ be the set of all 
$[\om]\in H_{\textrm{\tiny deRham}}^1(M,\mR)$ satisfying 
  $$\int_X \om \in r\mZ\quad\mbox{for any closed 1-chain }X.$$
Generalizing our definition in section~\ref{systolenormsection}
we define
\begin{eqnarray*}
  P:\homgitter{{1\over 2}}&\to &\Hom_\mZ(H_1(M,\mZ),\mZ_2)=H^1(M,\mZ_2)\\
  {}[\om] & \mapsto & \left([X] \mapsto \exp (2\pi i \int_X \om)\right).
\end{eqnarray*}
The kernel of $P$ is $\homgitter{}$.
We now define the \emph{stable norm} for elements of $\chi$ of $H^1(M,\mZ_2)$
  $$\lnorm{\chi}\infty := \inf \left\{ \lnorm{\om}\infty \,|\, P([\om])=\chi\right\}.$$
In general $P$ is not surjective, hence this norm takes values in $[0,\infty]$.
The elements in the image of $P$ are called {\em realizable by a dif\-fe\-ren\-tiab\-le form}.
A homomorphism $\chi\in \Hom_\mZ(H_1(M,\mZ),\mZ_2)$ is realizable by a dif\-fe\-ren\-tiab\-le form
if and only if $\chi$ vanishes on the torsion subgroup of $H_1(M,\mZ)$.

\begin{definition}
Two families $(\la_i|i\in\mZ)$ and $(\la_i'|i\in\mZ)$ of real numbers are said to be \emph{$\de$-close}
if there is a bijective map $h:\mZ\to \mZ$ with the property
  $$|\la_{h(i)}-\la_i'|\leq \de.$$
\end{definition}

\begin{proposition}\label{beldimprop}
Assume that $(M,g)$ carries two spin structures whose difference $\chi$ is realizable as a differential form.
Then the spectra of $D$ for the two spin structures are $2\pi\lnorm{\chi}\infty$-close.
\end{proposition}

\proof{} We modify a technique used by Friedrich \cite{friedrich:84} for calculating 
the spectrum of the Dirac operator on a flat torus.

Let us assume that the difference $\chi$ of the spin structures is realizable as a 
differentiable form.
We take $\om\in\homgitter{(1/2)}$ with $P([\om])=\chi$ and 
$\lnorm{\om}\infty\leq \lnorm{\chi}\infty+\ep$ for a small number $\ep>0$.
Then there is complex line bundle $L_\om$ on $M$ which is trivialized by a section
$\tau$ and a connection $\na$ on $L_\om$ such that
  $$\na_X\tau=2\pi i \om(X)\tau.$$
The holonomy of the bundle ($L_\om,\na)$ is $\chi$. 
Therefore the square of $(L_\om,\na)$ admits 
a parallel trivialization. 
Let $L_\om$ carry the hermitian metric characterized by $|\tau|\equiv 1$.

Denote by $\Si M$ and $\Si'M$ the spinor bundles to the two spin structures. Then
  $$\Si'M\cong\Si M\otimes L_\om$$
where the isomorphism preserves the connection, the hermitian metric and the Clifford 
multiplication. 
Now we define
\begin{eqnarray*} 
H:\Ga(\Si M)& \to &\Ga(\Si'M)\\
\Psi &\mapsto& \Psi\otimes \tau
\end{eqnarray*}
The Dirac operators $D$ and $D'$ for the two spin structures then satisfy
  $$D'\Psi =H\circ D \circ H^{-1} \Psi + 2\pi i \om \cdot\Psi$$ 
where $\cdot$ denotes the Clifford multiplication. Multiplication by $2\pi i \om$ is
a bounded operator on the space of $L^2$-sections of $\Si' M$. 
Its operator norm is $2\pi\lnorm{\om}\infty$. 
The following well-known lemma completes the proof.
\qed
\begin{lemma}
Let $D$ and $D'$ be two self-adjoint densely defined endomorphisms
of a complex separable Hilbert space.  
We assume that the spectra of $D$ and $D'$ are discrete with finite multiplicities.
Suppose that $D-D'$ is a bounded operator of 
operator norm $K$. Then the spectra of $D$ and $D'$ (with multiplicities) are $K$-close. 
\end{lemma}
The lemma is well-known in perturbation theory.  For example it can be deduced 
from considerations in \cite{kato:66}. 
The eigenvalues $(\la_i(t)\,|\,i\in\mZ)$ of 
  $$A_t:=(1-t) D + t D', \quad t\in [0,1]$$
can be numbered such that $\la_i(t)$ is a Lipschitz function in $t$ with Lipschitz constant $K$.
From this observation the lemma is evident.

\subsection*{Acknowledgement}
I want to thank Christian B\"ar for many interesting and stimulating discussions about the subject.
Christian B\"ar was the adviser of my PhD thesis and
many results in this article are based on the work in my PhD thesis.
The article was finally completed while I enjoyed the 
hospitality of the Graduate Center of CUNY, New York. 
It is a pleasure to thank Jozef Dodziuk for many helpful comments.
\bibliographystyle{alpha}

\newpage
\vspace{1cm}
Author's address:
\vspace{5mm}
\parskip0ex

\leavevmode
\komment{
\vtop{\hsize=8cm\noindent
{\obeylines
\emph{Until August 2000:}
\smallskip
CUNY Graduate Center
PhD Program in Mathematics
365, Fifth Avenue
New York, NY 10016
USA
}}
\vtop{\hsize=8cm\noindent
{\obeylines
\emph{From September 2000:}
\smallskip
Mathematisches Seminar
Universit\"at Hamburg
Bundesstra\ss{}e 55
20\,146 Hamburg
Germany
}}
}
{\obeylines
Bernd Ammann
Mathematisches Seminar
Universit\"at Hamburg
Bundesstra\ss{}e 55
20\,146 Hamburg
Germany
}\vspace{0.5cm}

E-Mail:
{\tt ammann@math.uni-hamburg.de}

WWW:
{\tt http://www.math.uni-hamburg.de/home/ammann}

\end{document}